\documentclass{gtart_h}  

\def\ifplaintex{\expandafter\ifx\csname documentclass\endcsname\relax}

\ifplaintex 
\else
\fi

\expandafter\ifx\csname beginpicture\endcsname\relax
\expandafter\ifx\csname documentclass\endcsname\relax
\input pictex \else
\input prepictex \input pictex \input postpictex \fi\fi

\def\gt{{\mathsurround=0pt\it $\cal G\mskip-2mu$eometry \&\ 
$\cal T\!\!$opology}}        

\def\gtp{{\mathsurround=0pt\it $\cal G\mskip-2mu$eometry \&\ 
$\cal T\!\!$opology $\cal P\!$ublications}}  

\def\lognumber#1{\def\thelognumber{#1}}
\def\volumenumber#1{\def\thevolumenumber{#1}}
\def\papernumber#1{\def\thepapernumber{#1}}
\def\volumeyear#1{\def\thevolumeyear{#1}}

\def\pagenumbers#1#2{\def\startpage{#1}\def\finishpage{#2}}
\def\published#1{\def\publishdate{#1}}
\def\proposed#1{\def\theproposer{#1}}
\def\seconded#1{\def\theseconders{#1}}
\def\received#1{\def\receiveddate{#1}}
\def\revised#1{\def\reviseddate{#1}}
\def\accepted#1{\def\accepteddate{#1}}
\def\asciititle#1{\def\theasciititle{#1}}
\def\covertitle#1{\def\thecovertitle{#1}}

\long\def\asciiabstract#1{\long\def\theasciiabstract{#1}}
\def\asciikeywords#1{\def\theasciikeywords{#1}}

\let\\\par\let\thelognumber\relax
\let\thevolumenumber\relax\let\thepapernumber\relax
\let\thevolumeyear\relax\let\thesamplenumber\relax\let\startpage\relax
\let\finishpage\relax\let\publishdate\relax\let\receiveddate\relax
\let\reviseddate\relax\let\accepteddate\relax\let\theasciititle\relax
\let\thecovertitle\relax\let\theasciiauthors\relax
\let\theasciiabstract\relax\let\theasciikeywords\relax
\let\theasciiemail\relax\let\theshortauthors\relax\let\theshorttitle\relax

\long\def\maketitlep{   

\count0=\startpage

\gt\hfill      
\beginpicture
\setcoordinatesystem units <0.33truein, 0.33truein> point at 2.2 0.9
\setplotsymbol ({$\cal G$})
\plotsymbolspacing=9truept
\circulararc 315 degrees from 0 1 center at 0 0
\setplotsymbol ({$\cal T$})
\circulararc 315 degrees from 1 -1 center at 1 0
\endpicture
\break
{\small\ifx\thesamplenumber\relax 
Volume \else Sample
\fi\thevolumenumber\ (\thevolumeyear)
\startpage--\finishpage\nl
Published: \publishdate}
\vglue 0.5truein plus 0.4fil minus 0.1truein

{\parskip=0pt\leftskip 0pt plus 1fil\def\\{\par\smallskip}{\ifplaintex\large
\else\Large\fi\bf\thetitle}\par\medskip}   

\vglue 0pt plus 0.1fil 

{\parskip=0pt\leftskip 0pt plus 1fil\def\\{\par}{\sc\theauthors}
\par\medskip}

\vglue 0pt plus 0.1fil 

{\small\parskip=0pt\let\newline\\
{\leftskip 0pt plus 1fil\def\\{\par}{\sl\theaddress}\par}
\expandafter\ifx\theemail\relax    
\relax\else\vglue 5pt plus 0.02fil minus 2pt\def\\{\stdspace{\rm 
and}\stdspace} 
\cl{Email:\stdspace\tt\theemail}\fi
\ifx\theurl\relax                  
\relax\else\vglue 5pt plus 0.02fil minus 2pt\def\\{\stdspace{\rm 
and}\stdspace}
\cl{URL:\stdspace\tt\theurl}\fi\par}

\vglue 7pt plus 0.3fil minus 3pt

{\bf Abstract}
\vglue 5pt plus 0.1fil minus 2pt

\theabstract

\vglue 7pt plus 0.3fil minus 3pt

{\bf AMS Classification numbers}\quad Primary:\quad \theprimaryclass

Secondary:\quad \thesecondaryclass

\vglue 5pt plus 0.3fil minus 2pt

{\bf Keywords:}\quad \thekeywords

\vglue 10pt plus 0.5fil minus 5pt

{\small  Proposed: \theproposer\hfill Received: \receiveddate\nl
Seconded: \theseconders\hfill 
\ifx\reviseddate\relax                         
Accepted: \accepteddate                        
\else
Revised: \reviseddate                          
\fi}
\eject
}       

\let\maketitlepage\maketitlep
\let\maketitle\maketitlepage

\font\phead=cmsl9 scaled 950
\font=cmsl9 scaled 1050
\font\pnum=cmbx10 scaled 913
\font\lnum=cmbx10 
\font\pfoot=cmsl9 scaled 950
\font=cmsl9 scaled 1050
\ifplaintex
\headline{\vbox to 0pt{\vskip -4.5mm\line{\small\phead\ifnum
\count0=\startpage ISSN 1364-0380 (on line)
1465-3060 (printed) \hfill {\pnum\folio}\else\ifodd\count0\def\\{ }%
\ifx\theshorttitle\relax\thetitle\else\theshorttitle\fi\hfill{\pnum\folio}
\else\def\\{ and }{\pnum\folio}\hfill\ifx\theshortauthors\relax\theauthors
\else\theshortauthors\fi\fi\fi}\vss}}
\footline{\vbox to 0pt{\vglue 0mm\line{\small\pfoot\ifnum\count0=\startpage
\copyright\ \gtp\hfill\else
\gt, Volume \thevolumenumber\ (\thevolumeyear)\hfill\fi}\vss
}}
\else
\makeatletter
\def\@oddhead{{\small\lhead\ifnum\count0=\startpage ISSN 1364-0380 (on line)
1465-3060 (printed) \hfill {\lnum\number\count0}\else\ifodd\count0
\def\\{ }\ifx\theshorttitle\relax \thetitle \else\theshorttitle\fi\hfill
{\lnum\number\count0}\else\def\\{ and }{\lnum\number\count0}
\hfill\ifx\theshortauthors\relax 
\theauthors\else\theshortauthors\fi\fi\fi}}\def\@evenhead{\@oddhead}
\def\@oddfoot{\small\lfoot\ifnum\count0=\startpage\copyright\ \gtp\hfill\else
\gt, Volume \thevolumenumber\ (\thevolumeyear)\hfill\fi}
\def\@evenfoot{\@oddfoot}
\makeatother
\fi

\newwrite\gtoutfile
\long\gdef\makeheadfile{  
{\def\\{, }\def\s{ }
\immediate\openout\gtoutfile head.xxx
\immediate\write\gtoutfile{Proxy-for: \ifx\theasciiauthors\relax
\theauthors\else\theasciiauthors\fi\s<\ifx\theasciiemail\relax\theemail\else\theasciiemail\fi>}
\immediate\write\gtoutfile{\noexpand\\}
\immediate\write\gtoutfile{Authors: \ifx\theasciiauthors\relax
\theauthors\else\theasciiauthors\fi}
{\def\\{ }\immediate\write\gtoutfile{Title: \ifx\theasciititle\relax
\thetitle\else\theasciititle\fi}}
\immediate\write\gtoutfile{Subj-class: GT or SG or MG etc}
\immediate\write\gtoutfile{MSC-class: \theprimaryclass\ifx\thesecondaryclass\relax\else, \thesecondaryclass\fi}
\immediate\write\gtoutfile{Journal-ref: Geom. Topol. \thevolumenumber
(\thevolumeyear) \startpage-\finishpage}
\immediate\write\gtoutfile{Comments: Published by Geometry and Topology at}
\immediate\write\gtoutfile{\s\s http://www.maths.warwick.ac.uk/gt/GTVol\thevolumenumber/paper\thepapernumber.abs.html}
\immediate\write\gtoutfile{\noexpand\\}
\immediate\write\gtoutfile{}
\ifx\theasciiabstract\relax
\immediate\write\gtoutfile{\theabstract}\else
\immediate\write\gtoutfile{\theasciiabstract}\fi
\immediate\write\gtoutfile{}
\immediate\write\gtoutfile{\noexpand\\}
\immediate\write\gtoutfile{}
\immediate\closeout\gtoutfile}}  

\def\maketitlepage{\maketitlep\makeheadfile}
\let\maketitle\maketitlepage

\lognumber{423}
\received{16 February 2004}
\volumenumber{8}\papernumber{36}\volumeyear{2004}
\pagenumbers{1301}{1359}
\revised{17 August 2004}
\published{19 October 2004}
\accepted{11 October 2004}
\proposed{Benson Farb}
\seconded{Walter Neumann, Joan Birman}

\usepackage{amsfonts, epsfig, amscd, amsmath}
\newtheorem{theorem}{Theorem}[section]
\newtheorem{lemma}[theorem]{Lemma}
\newtheorem{corollary}[theorem]{Corollary}

\newtheorem{question}[theorem]{Question}

\newtheorem{proposition}[theorem]{Proposition}

\def\Mod{\mbox{\rm{Mod}}}
\def\Homeo{\mbox{\rm{Homeo}}}

\def\Stab{\mbox{\rm{Stab}}}
\def\SL{\mbox{\rm{SL}}}
\def\PSL{\mbox{\rm{PSL}}}
\def\SO{\mbox{\rm{SO}}}
\def\DAf{\mbox{\rm{DAf}}}

\def\On{\mbox{\rm{O}}}
\def\length{\mbox{\rm{length}}}
\def\Aut{\mbox{\rm{Aut}}}
\def\Tr{\mbox{\rm{Tr}}}
\begin{document}
%
%
\title[On groups generated by two positive multi-twists]{On groups generated by two positive multi-twists:\\Teichm\"uller curves and Lehmer's number}

\covertitle{On groups generated by two positive multi-twists:\\Teichm\noexpand\"uller curves and Lehmer's number}

\asciititle{On groups generated by two positive multi-twists:\\Teichmueller curves and Lehmer's number}

\author{Christopher J Leininger}
\address{Department of Mathematics, Columbia University\\2990 Broadway MC 4448, New York, NY 10027, USA}
\email{clein@math.columbia.edu}
\begin{abstract}
From a simple observation about a construction of Thurston, we derive
several interesting facts about subgroups of the mapping class group
generated by two positive multi-twists.  In particular, we identify
all configurations of curves for which the corresponding groups fail
to be free, and show that a subset of these determine the same set
of Teichm\"uller curves as the non-obtuse lattice triangles which
were classified by Kenyon, Smillie, and Puchta.  We also identify a
pseudo-Anosov automorphism whose dilatation is Lehmer's number, and show
that this is minimal for the groups under consideration.  In addition,
we describe a connection to work of McMullen on Coxeter groups and
related work of Hironaka on a construction of an interesting class of
fibered links.
\end{abstract}
\asciiabstract{%
From a simple observation about a construction of Thurston, we derive
several interesting facts about subgroups of the mapping class group
generated by two positive multi-twists.  In particular, we identify
all configurations of curves for which the corresponding groups fail
to be free, and show that a subset of these determine the same set
of Teichmueller curves as the non-obtuse lattice triangles which
were classified by Kenyon, Smillie, and Puchta.  We also identify a
pseudo-Anosov automorphism whose dilatation is Lehmer's number, and show
that this is minimal for the groups under consideration.  In addition,
we describe a connection to work of McMullen on Coxeter groups and
related work of Hironaka on a construction of an interesting class of
fibered links.}
\asciikeywords{Coxeter, Dehn twist, Lehmer, pseudo-Anosov, mapping class group,
Teichmueller}
%
%
%
\maketitlepage    
%
%
%
\section{Introduction} \label{introsection}

Let $S$ be a connected finite type oriented surface.  In $\Mod(S)$,
the mapping class group of $S$, a particularly tractable class of
elements (or {\em automorphisms}) are the {\em positive multi-twists}.
These are products of positive Dehn twists about disjoint essential simple
closed curves.  For a given positive multi-twist, the union of these
simple closed curves is a closed essential $1$--manifold, and the set
of positive multi-twists is in a one-to-one correspondence with ${\cal
S}'(S)$, the set of isotopy classes of essential $1$--manifolds on $S$.
Given $A \in {\cal S}'(S)$, we let $T_{A}$ denote the positive multi-twist
which is the product of positive Dehn twists about the components of $A$.

This paper is concerned with subgroups of $\Mod(S)$ generated by two
positive multi-twists and is based on a construction of Thurston
\cite{T} (see also Long \cite{L} and Veech \cite{V1}).  When $A \cup B$
{\em fills} the surface (that is, every essential curve intersects $A$ or
$B$) Thurston constructs a certain type of Euclidean cone metric, which
we refer to as a {\em flat structure}, for which $\langle T_{A},T_{B}
\rangle$ acts by affine homeomorphisms.  The derivative of this action
defines a discrete homomorphism $\DAf\co \langle T_{A},T_{B} \rangle
\rightarrow \PSL_{2}({\mathbb R})$ with finite kernel.  This homomorphism
is determined by a single number, $\mu(A \cup B)$, depending on the
geometric intersection numbers of the components of $A$ with those of $B$.

\begin{figure}[ht!]
\begin{center}
\psfig{file=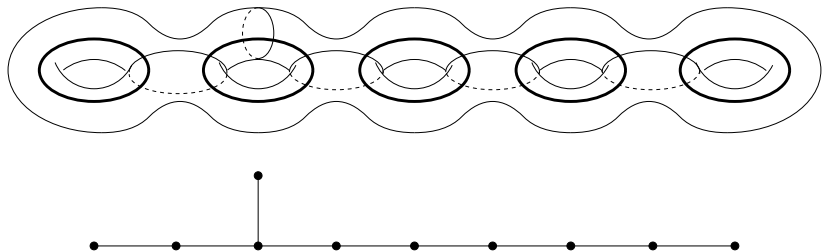,height=1.3truein}
\caption{$1$--manifolds $A_{L}$ and $B_{L}$ with configuration graph
${\cal G}(A_{L} \cup B_{L})={\cal E}h_{10}$}
\label{lehmerscurves}
\end{center}
\end{figure}

The novelty in our approach to studying these groups is Proposition
\ref{whatismu} in which we show that $\mu(A \cup B)$ is the spectral
radius of the {\em configuration graph}, ${\cal G}(A \cup B)$.
This graph has a vertex for each component of $A$ and of $B$ and an
edge for every point of intersection between corresponding components
(see Figure \ref{lehmerscurves}).

This observation, along with some elementary hyperbolic geometry and
well-known graph theoretic results, has many interesting consequences.

\subsection{Freeness} \label{introfreenesssection}

\begin{figure}[ht!]
\begin{center}
\psfig{file=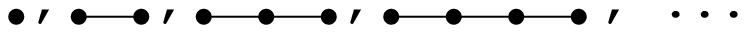,height=.2truein}
\caption{${\cal A}_{c}$, $c \geq 1$}
\label{straightchains}
\end{center}

\begin{center}
\psfig{file=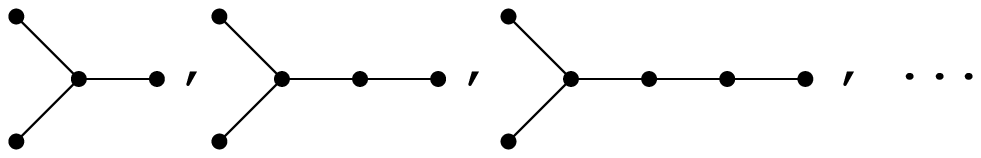,height=.5truein}
\caption{${\cal D}_{c}$, $c \geq 4$}
\label{hangers1}
\end{center}

\begin{center}
\psfig{file=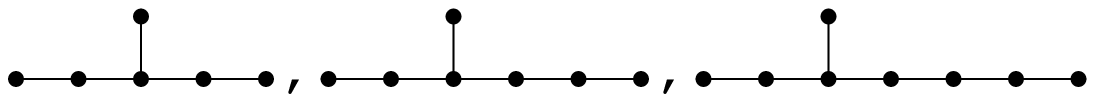,height=.35truein}
\caption{${\cal E}_{6}$, ${\cal E}_{7}$, and ${\cal E}_{8}$}
\label{exceptionals1}
\end{center}
\end{figure}

The graphs of type ${\cal A}_{c}$ ($c \geq 1$), ${\cal D}_{c}$ ($c
\geq 4$), ${\cal E}_{6}$, ${\cal E}_{7}$, and ${\cal E}_{8}$ play an
important role in our work, and we refer to them as {\em recessive}
graphs (see Figures \ref{straightchains}--\ref{exceptionals1}).
Any graph which is not recessive will be called {\em dominant}.\\

\medskip{\bf Theorem \ref{main1}}\qua
{\sl $\langle T_{A},T_{B} \rangle \cong {\mathbb F}_{2}$ if and only if ${\cal
G}(A \cup B)$ contains a dominant component.}\medskip

This theorem was inspired by the work of Hamidi-Tehrani in \cite{H} where
sufficient conditions for $\langle T_{A},T_{B} \rangle \cong {\mathbb
F}_{2}$ in terms of intersection numbers of components of $A$ with those
of $B$ are given.  In \cite{T}, Thurston remarks (without proof) that
a necessary and sufficient condition for this group to be free is that
$\mu(A \cup B) \geq 2$, and this is the basis for Theorem \ref{main1}.

\subsection{Teichm\"uller curves} \label{introteichcurvesection}

The proof of Theorem \ref{main1} reduces to the case that $A \cup
B$ fills $S$ (see Proposition \ref{reductionproposition} and Section
\ref{appendixreduction}), and we assume this to be the case for the
remainder of Section \ref{introsection}.

The flat structure on $S$ determines a quadratic differential and thus a
Teichm\"uller disk (see Section \ref{teichmullerspacesection}).  The group
$\langle T_{A},T_{B} \rangle$ stabilizes this disk, though in general
this group has infinite index in the full stabilizer.  The quotient of
a Teichm\"uller disk by its stabilizer is called a Teichm\"uller curve
when it has finite area, that is when the stabilizer is a lattice.
In this case, the Teichm\"uller curve isometrically immerses into the
moduli space of $S$, and we say that the Teichm\"uller disk covers a
Teichm\"uller curve.

In Zemljakov and Katok \cite{ZK} it is shown how to associate a Riemann
surface and quadratic differential to a rational polygon in such a way
that billiards trajectories in the polygon correspond to geodesics for
the associated flat structure (see Section \ref{billiardsection}).
A theorem of Veech \cite{V1} implies that the billiards in a polygon
have optimal dynamical properties if the corresponding Teichm\"uller
disk covers a Teichm\"uller curve, in which case the polygon is called
a {\em lattice polygon}.  In particular, understanding and classifying
Teichm\"uller curves and lattice polygons is an interesting problem which
has received much attention (see eg Veech \cite{V1}, \cite{V2}, Harvey
\cite{Har}, Gutkin and Judge \cite{GJ}, Kenyon and Smillie \cite{KS},
Puchta \cite{P}, McMullen \cite{Mc3}, and Calta \cite{Ca}).  Our second
main theorem provides a complete classification for a certain class of
Teichm\"uller curves.

\medskip
{\bf Theorem \ref{main2}}\qua
{\sl The Teichm\"uller curves for which the associated stabilizers contain
a group generated by two positive multi-twists with finite index are
precisely those defined by $A \cup B$ filling $S$, where ${\cal G}(A
\cup B)$ is critical or recessive.}\medskip

\begin{figure}[h]
\begin{center}
\psfig{file=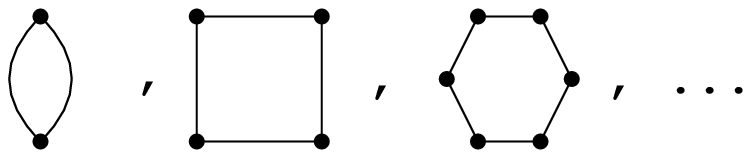,height=.5truein}
\caption{${\cal P}_{2c}$, $c \geq 1$}
\label{cyclicchains}
\end{center}

\begin{center}
\psfig{file=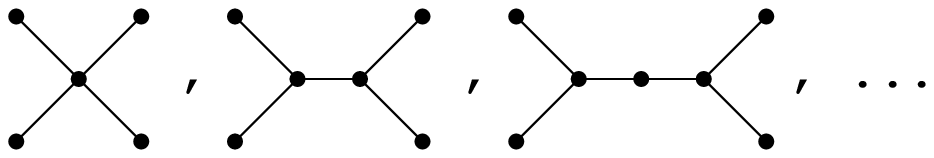,height=.5truein}
\caption{${\cal Q}_{c}$, $c \geq 5$}
\label{hangers2}
\end{center}

\begin{center}
\psfig{file=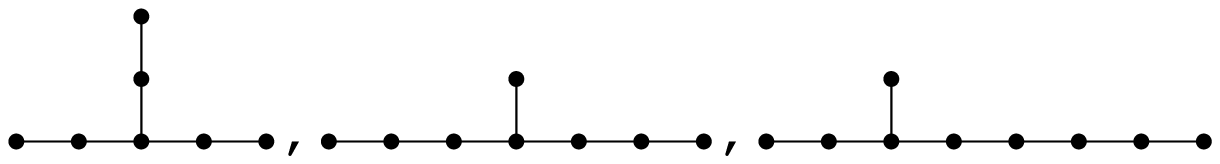,height=.5truein}
\caption{${\cal R}_{7}$, ${\cal R}_{8}$, and ${\cal R}_{9}$}
\label{exceptionals2}
\end{center}
\end{figure}

The critical graphs are those of type ${\cal P}_{2c}$ ($c \geq 1$), ${\cal
Q}_{c}$ ($c \geq 5$), ${\cal R}_{7}$, ${\cal R}_{8}$, and ${\cal R}_{9}$
(see Figures \ref{cyclicchains} -- \ref{exceptionals2}).

A classification of right and acute lattice triangles was initiated by
Kenyon and Smillie in \cite{KS} and completed by Puchta in \cite{Pu}
(see Theorem \ref{ksptheorem}).  Using this classification, we prove
the following theorem.

\medskip{\bf Theorem \ref{main3}}\qua
{\sl The Teichm\"uller curves determined by the right and acute lattice
triangles have associated stabilizers containing a finite index subgroup
of the form $\langle T_{A},T_{B} \rangle$ with ${\cal G}(A \cup B)$
recessive.}\medskip

Moreover, all the flat structures associated to ${\cal G}(A \cup B)$
of type ${\cal A}_{c}$ and ${\cal D}_{c}$ are affine equivalent to
structures which can be tiled by one or two regular Euclidean polygons
(see Section \ref{whichgroupsection}).  These were all studied by Veech
in \cite{V1} and \cite{V2} and by Earle and Gardiner in \cite{EG}.
The remaining three cases where ${\cal G}(A \cup B)$ has type ${\cal
E}_{6}$, ${\cal E}_{7}$, and ${\cal E}_{8}$ correspond to the exceptional
triangles mentioned above and were studied by Veech in \cite{V1}, Vorobets
in \cite{Vor}, and Kenyon and Smillie in \cite{KS}.  In particular,
in those cases that $\langle T_{A},T_{B} \rangle$ fails to be free,
the aforementioned references provide a description of these groups.

\medskip{\bf Theorem \ref{main4}}\qua
{\sl If ${\cal G}(A \cup B)$ is recessive, then $\DAf$ maps $\langle
T_{A},T_{B} \rangle$ onto a Fuchsian triangle group with finite central
kernel of order at most $2$.
The signature of the triangle group is described by the following table.
\def\strutt{\vrule width 0pt height 15pt depth 5pt}
$$\begin{array} {cccc}
\strutt\mbox{configuration graph} & \mbox{signature} & \mbox{configuration graph} & \mbox{signature} \\
\hline
\strutt{\cal D}_{c}, \, c \geq 4 & (c-1,\infty,\infty) & {\cal E}_{6} & (6,\infty,\infty)\\
\strutt{\cal A}_{2c+1}, \, c \geq 1 & (c+1,\infty,\infty) & {\cal E}_{7} & (9,\infty,\infty)\\
\strutt{\cal A}_{2c}, \, c \geq 1 & (2,2c+1,\infty) & {\cal E}_{8} & (15,\infty,\infty)\\ \end{array}$$}

\subsection{Lehmer's number and Coxeter groups} \label{introlehmercoxsection}

The original purpose of Thurston's construction was not to study
groups generated by two positive multi-twists, but rather to construct
explicit examples of {\em pseudo-Anosov} automorphisms.  Indeed, in
these groups, pseudo-Anosov automorphisms are generic (see Proposition
\ref{elementtypes}).  Associated to a pseudo-Anosov automorphism, $\phi$,
is an algebraic integer $\lambda(\phi) > 1$ called the {\em dilatation}
reflecting certain dynamical properties (see Section \ref{pseudosection}).

\medskip{\bf Theorem \ref{main5}}\qua
{\sl For any surface $S$, any $A, B \in {\cal S}'(S)$, and any
pseudo-Anosov element
$$\phi \in \langle T_{A},T_{B} \rangle <  \Mod(S)$$
we have $\lambda(\phi) \geq \lambda_{L} \approx 1.1762808$.
Moreover, $\lambda(\phi) = \lambda_{L}$ precisely when $S$ has genus
$5$ (with at most one marked point), $\{ A,B \} = \{ A_{L}, B_{L} \}$
as in Figure \ref{lehmerscurves} (up to homeomorphism), and $\phi$
is conjugate to $(T_{A}T_{B})^{\pm 1}$.}\medskip

Here $\lambda_{L}$ is {\em Lehmer's number} which is the largest real
root of {\em Lehmer's polynomial}:
\begin{equation} \label{lehmerseqn}
x^{10} + x^{9} - x^{7} - x^{6} - x^{5} - x^{4} - x^{3} + x + 1
\end{equation}
$\lambda_{L}$ was discovered by Lehmer in 1933 \cite{Leh} and is the
smallest known Salem number and Mahler measure of an integral polynomial
(see Section \ref{lehmersapplication}).

Given any ${\cal G} = {\cal G}(A \cup B)$, we can view this as a Coxeter
graph, and we let ${\frak C}({\cal G})$ and ${\frak A}({\cal G})$ denote
the associated {\em Coxeter group} and {\em Artin group}, respectively,
with $\pi_{ac}\co{\frak A}({\cal G}) \rightarrow {\frak C}({\cal G})$
the canonical epimorphism.  We let $\Theta$ denote the geometric
action of ${\frak C}({\cal G})$ on $({\mathbb R}^{K},\Pi_{{\cal G}})$
(see Section \ref{coxartsection}).

Theorem \ref{main5} is strikingly similar to the main theorem of
\cite{Mc} (see Theorem \ref{mctheorem}) in which McMullen shows that the
minimal spectral radius of any essential element in a Coxeter group with
respect to $\Theta$ is either $1$ or else bounded below by $\lambda_{L}$.
Moreover, $\lambda_{L}$ is achieved precisely when the associated Coxeter
graph is ${\cal E}h_{10}$.

We say ${\cal G}$ has {\em small type} if there are no multiple edges
between vertices.  It is well known that there is a homomorphism
$$\Psi\co{\frak A}({\cal G}(A \cup B)) \rightarrow \Mod(S)$$ sending
the standard generators to the corresponding Dehn twists in the $A$
and $B$ curves, when ${\cal G}(A \cup B)$ has small type (see Section
\ref{artmapsection}).  This provides the first link with the groups under
consideration.  The following describes the connection with McMullen's
Theorem \ref{mctheorem}.

\medskip{\bf Theorem \ref{main6}}\qua
{\sl Let ${\cal G}(A \cup B)$ be non-critical dominant with small type.
Then $\sigma_{A}\sigma_{B}$ is sent by $\Psi$ to a pseudo-Anosov with
dilatation equal to the spectral radius of its image under $\Theta
\circ \pi_{ac}$.  Moreover, among all essential elements in $\langle
\sigma_{A},\sigma_{B} \rangle$, $\sigma_{A}\sigma_{B}$ minimizes both
dilatation as well as spectral radius for the respective homomorphisms.}
\medskip

In this theorem, $\sigma_{A}\sigma_{B}$ is the {\em bicolored Coxeter
element} (see Section \ref{coxbasics}).  Inspired by the work of
Hironaka in \cite{Hir1} (see Theorem \ref{hirtheorem}), we find that
under some additional hypothesis, (part of) the action of $T_{A}T_{B}$
on $H_{1}(S;{\mathbb R})$ is almost semi-conjugate to the geometric
action of this Coxeter element.

\medskip{\bf Theorem \ref{main7}}\qua
{\sl Let ${\cal G}(A \cup B)$ have small type and suppose that $A$ and $B$
can be oriented so that all intersections of $A$ with $B$ are positive.
Then there exists a homomorphism
$$\eta\co{\mathbb R}^{K} \rightarrow H_{1}(S;{\mathbb R})$$
such that
$$(T_{A}T_{B})_{*} \circ \eta = - \eta \circ \Theta(\sigma_{A}\sigma_{B}).$$
Moreover, $\Theta(\sigma_{A}\sigma_{B})|_{\ker(\eta)} = -I$ and $\eta$
preserves spectral radii.}\medskip

To say that $\eta$ preserves spectral radii, we simply mean that
$\Theta(\sigma_{A}\sigma_{B})$ and $(T_{A}T_{B})_{*}$ have the same
spectral radius (modulus of leading eigenvalue).  If we wish to relate
$T_{A}T_{B}$ to $\sigma_{A}\sigma_{B}$, this theorem is likely the best
we can do.  For, unlike Hironaka's Theorem \ref{hirtheorem}, there is
no relation between $\dim(H_{1}(S;{\mathbb R}))$ and $K$.

\medskip
{\bf Remarks}

(1)\qua Theorem \ref{main1} allows the possibility that $A \cup B$ does
not fill $S$, while we implicitly assume this for the other theorems.
We also note that Theorem \ref{main1} holds when $S$ has nonempty boundary
(see Section \ref{appendixreduction}). 

(2)\qua We caution the reader that the connection to Coxeter groups we
have described is only valid when the configuration graph has small type.
This is easily explained by the fact that the adjacency matrix for a
graph is the same as the {\em Coxeter adjacency matrix} only when the
graph has small type.

The paper is organized as follows.  Sections \ref{background}
and \ref{teichmullerspacesection} contain definitions and theorems
regarding surface topology, mapping class groups, and Teichm\"uller
space.  We recall the relevant facts concerning matrices and graphs in
Section \ref{pfsection}.  In Section \ref{singeuc} we give Thurston's
construction and prove Proposition \ref{whatismu} relating this to
the spectral radius of the configuration graph.  Next we discuss some
basics of Fuchsian groups and use them to prove Theorems \ref{main1} and
\ref{main5} in Section \ref{fuchsian}.  In Section \ref{whichgroupsection}
we discuss in more detail the groups corresponding to the critical and
recessive configurations and prove Theorems \ref{main2}, \ref{main3},
and \ref{main4}.  We then turn to Coxeter and Artin groups in Section
\ref{coxartsection}, describe the connection with groups generated by two
positive multi-twists, and prove Theorems \ref{main6} and \ref{main7}.
In Section \ref{applicationsection} we provide a few applications of
the theorems and indicate some interesting open questions.

We have also included two appendices.  The first, Section
\ref{appendixreduction}, reduces the proof of Theorem \ref{main1} to the
filling case, as well as extending it to the situation of surfaces with
boundary.  The second, Section \ref{pennerconst}, addresses a construction
of pseudo-Anosov automorphisms given by Penner which extends Thurston's
construction.  For completeness, we show that the lower bound given by
Theorem \ref{main5} holds for this class as well.

\eject

{\bf Acknowledgements}

I would like to thank Norbert A'Campo, Joan Birman, Hessam Hamidi-Tehrani,
Curt McMullen, Walter Neumann, and Alan Reid for helpful conversations and
some very useful and interesting references which not only simplified
the exposition, but which led to some of the connections to other
areas of mathematics.  Thanks also to Dan Margalit and the referees
for carefully reading earlier versions of this paper and providing
many helpful suggestions.  The author was partially supported by an
N.S.F. postdoctoral fellowship.

\section{Surface topology and mapping class groups} \label{background}

For more details on the material reviewed in this section, see Thurston
\cite{T}, Birman \cite{Bi}, Ivanov \cite{Iv}, and the lecture notes
\cite{FLP}.

\subsection{Surfaces and essential $1$--manifolds} \label{backsurfacesection}

Let $S = S_{g,p}$ be a smooth, compact, connected, oriented, genus--$g$
surface with $p$ marked points.  We will ignore the trivial cases,
and hence from this point on assume $S \neq S_{0,p}$ for $p \leq 3$.
Denote by $\dot{S}$ the surface $S$ with the $p$ marked points removed.

We denote the set of isotopy classes of essential simple closed curves on
$S$ by ${\cal S}(S)$.  That is, an element of ${\cal S}(S)$ is an isotopy
class of homotopically essential simple closed curves on $\dot{S}$ not
isotopic to a puncture of $\dot{S}$.  The geometric intersection number
for a pair of elements $a,b \in {\cal S}(S)$, denoted $i(a,b)$, is the
minimal number of transverse intersection points among all representatives
of $a$ and of $b$.

Let ${\cal S}'(S)$ denote the set of isotopy classes of essential,
closed $1$--manifolds embedded in $S$.  An element $A \in {\cal S}'(S)$
is an embedded $1$--submanifold of $\dot{S}$, for which every component
is homotopically essential in $\dot{S}$ and not isotopic to a puncture,
well-defined up to isotopy.  We will make no distinction between
$1$--manifolds and the isotopy classes they represent when convenient.
We refer to the components of $A$ as elements of ${\cal S}(S)$.  Whenever
we write $A = a_{1} \cup \cdots \cup a_{n}$ it will be assumed that $a_{i}
\in {\cal S}(S)$, for each $i = 1,\ldots,n$.

Note that an element of ${\cal S}'(S)$ is allowed to have several of
its components parallel (isotopic) to one another.

When considering elements $A,B \in {\cal S}'(S)$ as representative
$1$--manifolds of their isotopy classes, we will always assume that they
meet transversely and minimally.
When this is done, a component $a$ of $A$ and $b$ of $B$ meet in exactly
$i(a,b)$ points.
Consequently, the configuration graph ${\cal G}(A \cup B)$ depends only
on $A$ and $B$, and its components are in a one-to-one correspondence
with those of $A \cup B$, thought of as a subset of $S$.

We further note that ${\cal G}(A \cup B)$ is {\em bipartite}.
That is, the vertices may be colored by two colors (call them $A$ and $B$)
so that no two vertices of the same color are adjacent.

When $A \cup B$ fills $S$ (see Section \ref{introsection}) it follows
that the components of $S \setminus (A \cup B)$ are disks (each with
at most one marked point).  Note that ${\cal G}(A \cup B)$ is connected
when $A \cup B$ fills $S$.

\subsection{Uniqueness} \label{uniquesection}

Given a bipartite graph ${\cal G}$, there may be several different
pairs of $1$--manifolds having ${\cal G}$ as the configuration graph.
Indeed, ${\cal G}$ need not even determine the homeomorphism type of
the underlying surface.  However, there are instances in which one does
have uniqueness.

\begin{proposition} \label{recessiveuniqueness}
Suppose $A_{i},B_{i} \in {\cal S}'(S_{i})$, $A_{i} \cup B_{i}$ filling
$S_{i}$, $i = 1,2$, and ${\cal G} = {\cal G}(A_{1} \cup B_{1}) = {\cal
G}(A_{2} \cup B_{2})$.
If ${\cal G}$ is a tree with only one vertex of valence at most three
(in particular, if it is recessive), then there is a homeomorphism from
$S_{1}$ to $S_{2}$ taking $\{ A_{1},B_{1} \}$ to $\{ A_{2}, B_{2} \}$,
up to adding marked points.
\end{proposition}

\begin{figure}[ht!]\small
\begin{center}
\psfig{file=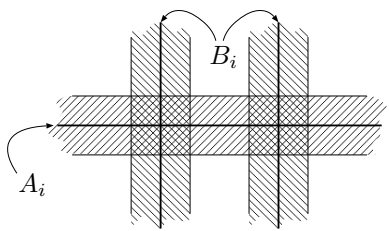,height=1.2truein}
\caption{Pieces of annuli intersecting in squares}
\label{intersecsquares}
\end{center}
  \setlength{\unitlength}{1in}
  \begin{picture}(0,0)(0,0)
    \put(1.7,.72){$A_{i}$}
    \put(2.7,1.39){$B_{i}$}
  \end{picture}
\end{figure}

{\bf Sketch of proof}\qua
Let ${\cal N}(A_{i} \cup B_{i})$ denote a regular neighborhood of $A_{i}
\cup B_{i}$ in $S_{i}$.
$S_{i}$ is obtained from ${\cal N}(A_{i} \cup B_{i})$ by adding disks
with zero or one marked point each, and there is just one way to do this,
up to homeomorphism.
So, to find a homeomorphism from $S_{1}$ to $S_{2}$, it suffices to find
a homeomorphism from ${\cal N}(A_{1} \cup B_{1})$ to ${\cal N}(A_{2}
\cup B_{2})$.
We view ${\cal N}(A_{i} \cup B_{i})$ as a union of annular neighborhoods
of the components of $A_{i}$ and $B_{i}$ with pairs of annuli intersecting
in squares if the corresponding curves intersect, and otherwise not at
all (see Figure \ref{intersecsquares}).
Thus, in each surface, we have an annulus associated to each vertex of
${\cal G}(A_{i} \cup B_{i})$ and a square of intersections of annuli
for each edge.

We first define a homeomorphism on the annulus corresponding to the
three-valent vertex, $v$, (if it exists, otherwise, we can start at any
vertex).  We do this so that the three (or fewer) squares corresponding
to the edges meeting $v$ are taken to the squares corresponding to the
same three edges.  This is possible because homeomorphisms of the circle
act transitively on ordered triples of points (and so homeomorphisms of
an annulus act transitively on ordered triples of disjoint squares).
Next we extend the homeomorphism over the annuli corresponding to the
vertices adjacent to $v$, again preserving the squares corresponding to
the edges meeting those vertices.  The homeomorphism is already defined
on one of these squares and now there is at most one other square (since
these vertices have valence at most $2$), so this is possible as well.
We may continue in this way extending over annuli corresponding to
adjacent vertices, preserving intersection squares.  At each stage,
we only encounter vertices with valence at most $2$, so this is always
possible.

Since there are only finitely many vertices, after finitely many steps
we obtain a homeomorphism from ${\cal N}(A_{1} \cup B_{1})$ to ${\cal
N}(A_{2} \cup B_{2})$, taking annuli to annuli.  The union of the cores
of these annuli is precisely $A_{i} \cup B_{i}$, so applying an isotopy
if necessary, we may assume that $A_{1} \cup B_{1}$ is taken to $A_{2}
\cup B_{2}$.
\endproof

\subsection{Automorphisms} \label{autosection}

We say that a homeomorphism $\phi\co S \rightarrow S$ is {\em allowable}
if it preserves the marked points.  We denote the group of allowable,
orientation preserving homeomorphisms of $S$ by $\Homeo_{+}(S)$ and
the identity component by $\Homeo_{0}(S)$.  The {\em mapping class
group} is defined to be the quotient group
$$\Mod(S) = \Homeo_{+}(S) / \Homeo_{0}(S).$$
An element of $\Mod(S)$ is referred to as an {\em automorphism} of $S$,
and by definition is a homeomorphism, well-defined up to isotopy.
When no confusion can arise, we will make no distinction between a
homeomorphism and the automorphism it determines.

\subsection{Multi-twists} \label{multitwistsection}

Given $a \in {\cal S}(S)$, a {\em positive Dehn twist} is the isotopy
class of a homeomorphism supported in an annular neighborhood of $a$
described as follows.
If we identify the annular neighborhood of $a$ with the annulus ${\mathbb
R}/\tau {\mathbb Z} \times [0,\sigma]$ by an orientation preserving
homeomorphism, then with respect to the obvious coordinates on this
annulus, $(t,s)$, the Dehn twist is given by
\begin{equation} \label{dehneqn}
(t,s) \mapsto \left(t + s \frac{\tau}{\sigma},s\right)
\end{equation}
We note that this makes the Dehn twist affine with respect to the natural
Euclidean metric on the annulus for any $\tau, \sigma > 0$.

Given $A = a_{1} \cup \cdots \cup a_{n} \in {\cal S}'(S)$, a {\em
multi-twist} along $A$ is the product
$$T_{a_{1}}^{\epsilon_{1}} \cdots T_{a_{n}}^{\epsilon_{n}}$$
where $\epsilon_{i} \in \{ \pm 1 \}$.  The {\em positive multi-twist}
along $A$, written $T_{A}$, is given by the above product where
$\epsilon_{i} = 1$ for each $i =1,\ldots,n$.  The map $A \mapsto T_{A}$
determines a bijection between ${\cal S}'(S)$ and the set of positive
multi-twists.\\

In the definition of $T_{A}$ the order of the product does not matter
since Dehn twists in disjoint curves obviously commute.
In fact, for $a,b \in {\cal S}(S)$, we have
\begin{equation} \label{braidrelationequation}
\begin{array}{ccc}
i(a,b) = 0 & \Rightarrow & T_{a}T_{b} = T_{b}T_{a}\\
i(a,b) = 1 & \Rightarrow & T_{a}T_{b}T_{a} = T_{b}T_{a}T_{b}\\
\end{array}
\end{equation}
The second equality is the well known {\em braid relation} and easily
follows from the calculation that $T_{a}T_{b}(a) = b$.

\medskip{\bf Remark}\qua
One often requires only that $\epsilon_{i} \in {\mathbb Z}_{+}$ in the
definition of a positive multi-twist.  However, we may replace a power of
a Dehn twist about a curve $a$ by a product of Dehn twist about several
isotopic copies of $a$, so there is no loss in generality in taking only
the powers $1$ in our definition.

\subsection{Pseudo-Anosov automorphisms} \label{pseudosection}

An automorphism $[\phi] \in \Mod(S)$ is called {\em pseudo-Anosov}
if there is a representative $\phi$ which leaves invariant a pair of
mutually transverse singular foliations with the following property.
These foliations admit transverse measures and $\phi$ multiplies one
measure by a factor $\lambda > 1$ and the other by $\lambda^{-1}$.
The number $\lambda = \lambda([\phi]) = \lambda(\phi) > 1$ is an algebraic
integer called the {\em dilatation} of $[\phi]$.

The dilatation has the following geometric description.  For any $a \in
{\cal S}(S)$ and any complete hyperbolic metric on $S$, the length of
the geodesic representative of $\phi^{n}(a)$ grows like $\lambda^{n}$.
That is, $\lambda^{-n} \length(\phi^{n}(a))$ converges to a nonzero
number.

\subsection{Reduction to the filling case} \label{reducingsection}

As was mentioned in the introduction, Theorem \ref{main1} is valid for
any surface and any pair of $1$--manifolds, but the proof reduces to
the filling case.  A proof of the following is given in the Appendix,
in Section \ref{appendixreduction}.

\medskip{\bf Proposition \ref{reductionproposition}}\qua
{\sl It suffices to prove Theorem \ref{main1} for $A \cup B$ filling $S$.}
\medskip

{\bf Convention}\qua For the remainder of this paper (excluding Section
\ref{appendixreduction}) we shall assume that every pair of essential
$1$--manifolds is filling.

\section{Teichm\"uller and moduli spaces} \label{teichmullerspacesection}

For more details on Teichm\"uller space and quadratic differentials see
Gardiner and Lakic \cite{Gard}, Masur \cite{Mas}, Masur and Tabachnikov
\cite{MT}, Earle and Gardiner \cite{EG}, and McMullen \cite{Mc3}.

Consider the space of complex structures on $S$, with orientation
compatible with the given orientation.  $\Homeo_{+}(S)$ acts on this
space, and the quotient is called the {\em moduli space} of $S$ and is
denoted ${\cal M}(S)$.  If we quotient by the action of the subgroup
$\Homeo_{0}(S)$ the resulting space is called the {\em Teichm\"uller
space} of $S$, and is denoted ${\cal T}(S)$.  ${\cal T}(S)$ is the
universal orbifold covering of ${\cal M}(S)$, with covering group
$\Mod(S)$.

Given $[J_{0}],[J_{1}] \in {\cal T}(S)$, the Teichm\"uller distance is
defined by
$$d([J_{0}],[J_{1}]) = \frac{1}{2}\inf_{f \simeq Id_{S}}
  \log\left(K\left(f\co (S,J_{0}) \rightarrow (S,J_{1})\right)\right),$$
where the infimum is taken over all quasi-conformal homeomorphisms $f$
isotopic to the identity, and $K(f)$ is the dilatation of $f$.
The action of $\Mod(S)$ is by isometries, and so the metric pushes down
to ${\cal M}(S)$.

\subsection{Quadratic differentials} \label{quaddiffsection}

Let $[J] \in {\cal T}(S)$, and consider the space ${\cal Q}(S,J)$
of integrable meromorphic quadratic differentials on $(S,J)$ which
are holomorphic on $\dot{S}$.  Any $q \in {\cal Q}(S,J)$ determines
a singular Euclidean metric $|q|$ on $S$ with cone-type singularities
having cone angles $k \pi$ for $k \in {\mathbb Z}_{\geq 3}$ at non-marked
points and $k \in {\mathbb Z}_{\geq 1}$ at marked points.  It also
defines a singular measured foliation ${\cal F}_{h}$, called the {\em
horizontal foliation}, whose leaves are geodesic with respect to $|q|$.
These leaves are precisely the injectively immersed $1$--manifolds $\gamma$
satisfying $q(\gamma'(t)) \geq 0$.  We refer to this structure as a {\em
flat structure}, and will also denote it by $q$.

When all the leaves of ${\cal F}_{h}$ are compact, the complement of the
singular leaves is a disjoint union of annuli.  In this situation, we say
that $q$ (or ${\cal F}_{h}$) determines an annular decomposition of $S$.

Suppose now we are given a flat structure $q$.  That is, we
have a singular Euclidean metric $|q|$ (having the above types of
singularities) and a singular foliation ${\cal F}_{h}$ with geodesic
leaves.  This defines a complex structure $J$ and quadratic differential
which can be described as follows.  The singular Euclidean metric is
given by an atlas of charts into ${\mathbb C}$ on the complement of the
singularities for which the transition functions are Euclidean isometries.
This defines a complex structure on the complement of the singularities
which then extends over this finite set.  Requiring that the leaves of
${\cal F}_{h}$ be sent to horizontal lines by our charts restricts our
transition functions to be of the form $z \mapsto \pm z + \xi$, for
some $\xi \in {\mathbb C}$.  We refer to such an atlas of charts as a
{\em preferred atlas} for $q$.  The form $dz^{2}$ is invariant under the
transition functions and pulls back to the desired quadratic differential.
The horizontal foliation is precisely ${\cal F}_{h}$.

We also obtain a locally defined orthonormal basis $e_{1},e_{2}$
for the tangent space to any non-singular point such that $e_{1}$
is tangent to ${\cal F}_{h}$.  Away from the singularities, this
basis is globally well-defined by this condition, up to sign (ie by
replacing $\{e_{1},e_{2}\}$ by $\{-e_{1},-e_{2}\}$).  The dual basis
$\{e^{1},e^{2}\}$ locally defines a holomorphic $1$--form $\omega = e^{1}
+ i e^{2}$.  Although $\omega$ is not in general globally well-defined,
its square is, and this is precisely the quadratic differential $q =
\omega^{2}$.  Note that $\omega$ is globally defined precisely when the
metric has no holonomy.

\subsection{Teichm\"uller disks and curves} \label{diskcurvesection}

Given $[J] \in {\cal T}(S)$, and $q \in {\cal Q}(S,J)$, there exists a map
$$\widetilde{f}\co\SL_{2}{\mathbb R} \rightarrow {\cal T}(S)$$
which sends $\gamma \in \SL_{2}{\mathbb R}$ to a point in ${\cal T}(S)$
obtained by deforming $[J]$ according to $\gamma$ as follows.  An element
$\gamma \in \SL_{2}{\mathbb R}$ defines a new atlas by composing each
chart in the preferred atlas with $\gamma$ (here we are identifying
${\mathbb C}$ with ${\mathbb R}^{2}$ and $\gamma$ is the obvious ${\mathbb
R}$--linear map).  The transition functions for the new atlas are again of
the form $z \mapsto \pm z + \xi$, and we obtain a new complex structure
$\gamma \cdot J$ and quadratic differential $\gamma \cdot q$.  We define
$\widetilde{f}(\gamma) = \gamma \cdot J$.

Note that $\SO(2)$ does not change the underlying complex structure,
and so $\tilde{f}$ factors through a map
$$f\co {\mathbb H}^{2} \cong \SO_{2} \setminus \SL_{2}({\mathbb R})
  \rightarrow {\cal T}(S).$$
After scaling the hyperbolic metric this is a holomorphic isometric
embedding and is called a {\em Teichm\"uller disk}.

Given a Teichm\"uller disk $f\co{\mathbb H}^{2} \rightarrow {\cal T}(S)$,
we have the stabilizer of $f({\mathbb H}^{2})$
$$\Stab(f({\mathbb H}^{2})) < \Mod(S).$$
Conjugating by $f$ we obtain a subgroup of $PSL_{2}({\mathbb R})$ which
we denote
$$\Stab(f) = f^{-1} \Stab(f({\mathbb H}^{2})) f.$$
Forming the quotient by $\Stab(f)$, $f$ then descends to a map
$$\widehat{f}\co {\mathbb H}^{2} / \Stab(f) \rightarrow {\cal M}(S).$$
When ${\mathbb H}^{2}/\Stab(f)$ has finite area, its image
$\widehat{f}({\mathbb H}^{2}/\Stab(f))$ is an algebraic curve totally
geodesically immersed in ${\cal M}(S)$ called a {\em Teichm\"uller curve}.

If $f$ is a Teichm\"uller disk defined by $q$, then every automorphism
of $\Stab(f({\mathbb H}^{2}))$ can be realized by an affine automorphism
with respect to the flat structure.
The derivative with respect to the basis $\{e_{1},e_{2} \}$ defines a
discrete representation
$$\DAf\co\Stab(f({\mathbb H}^{2})) \rightarrow \PSL_{2}({\mathbb R})$$
(this is into $\PSL_{2}({\mathbb R})$, rather than $\SL_{2}({\mathbb R})$
because the basis is only defined up to sign).

An element of the kernel of $\DAf$ leaves the complex structure and the
quadratic differential invariant.  It follows that such an element fixes
the Teichm\"uller disk pointwise.  Because the action on ${\cal T}(S)$
is properly discontinuous, the kernel of $\DAf$ is finite.

We collect these and other facts into the following theorem for ease of reference.
Parts of this theorem have appeared in several different locations (see eg the lecture notes \cite{FLP}, Thurston \cite{T}, Kra \cite{Kra}, Long \cite{L}, and Veech \cite{V1}).

\begin{theorem}[Thurston, Kra, Veech] \label{reptheorem}
Let $f\co{\mathbb H}^{2} \rightarrow {\cal T}(S)$ be a Teichm\"uller disk.
Then
$$\DAf\co\Stab(f({\mathbb H}^{2})) \rightarrow \PSL_{2}{\mathbb R}$$
is discrete, with finite kernel.
For $\phi \in \Stab(f({\mathbb H}^{2})) \setminus \{ 1 \}$ the following
is true:
\begin{enumerate}
\item if $\DAf(\phi)$ is elliptic or the identity, then $\phi$ has finite order,
\item if $\DAf(\phi)$ is parabolic, then $\phi$ is reducible and some
  power of $\phi$ is a positive multi-twist, and
\item if $\DAf(\phi)$ is hyperbolic, then $\phi$ is pseudo-Anosov and the
  dilatation is given by $\lambda(\phi) =
  \exp\left(\frac{1}{2}L(\DAf(\phi))\right)$ where $L(\DAf(\phi))$ is the
  translation length of $\DAf(\phi)$ on
${\mathbb H}^{2}$.
\end{enumerate}
\end{theorem}

There is a strong converse to part 3 of the theorem which is essentially Bers' description of pseudo-Anosov automorphisms.

\begin{theorem}[Bers] \label{bersuniquetheorem}
Given any pseudo-Anosov automorphism $\phi$, there is a unique
Teichm\"uller disk which it stabilizes.
\end{theorem}

The quotients ${\mathbb H}^{2}/\DAf(\Stab(f({\mathbb H}^{2})))$ and
${\mathbb H}^{2}/\Stab(f)$ are essentially the same (see eg \cite {EG}
or \cite{Mc3}).

\begin{proposition} \label{mcproposition}
${\mathbb H}^{2}/\Stab(f)$ is isometric to
${\mathbb H}^{2}/\DAf(\Stab(f({\mathbb H}^{2})))$.  In particular,
$\Stab(f)$ has finite co-area if and only if
$\DAf(\Stab(f({\mathbb H}^{2})))$ does.
\end{proposition}

\subsection{Homology representation} \label{homologysection}

As noted above, $q$ is the square of a holomorphic $1$--form $\omega
= e^{1} + i e^{2}$ if and only if the metric has no holonomy.
Now suppose $q = \omega^{2}$ and let $f\co{\mathbb H}^{2} \rightarrow
{\cal T}(S)$ denote the associated Teichm\"uller disk.  In this
case, the two-dimensional subspace $\langle e^{1},e^{2} \rangle
\subset H^{1}(S;{\mathbb R})$ is left invariant by the action of
$\Stab(f({\mathbb H}^{2}))$ since the ${\mathbb R}$--span of the vector
fields $\{e_{1},e_{2}\}$ is invariant.  Moreover, the action on $\langle
e^{1},e^{2} \rangle$ is dual to the action on $\langle e_{1},e_{2} \rangle
\cong {\mathbb R}^{2}$ given by $\DAf$.  In particular, the induced
action on $H^{1}(S;{\mathbb R})$ has a finite kernel.  Since finite
order automorphisms can never act trivially on (co)homology, we obtain
the following result.

\begin{proposition} \label{homologyreptheorem}
Suppose $q = \omega^{2}$ and $f\co{\mathbb H}^{2} \rightarrow {\cal T}(S)$
is the corresponding Teichm\"uller disk.  Then the action of
$\Stab(f({\mathbb H}^{2}))$ on homology is faithful.
\end{proposition}

\subsection{Example: Billiards} \label{billiardsection}

Let $P \subset {\mathbb R}^{2}$ be a compact {\em rational} polygon,
that is, the angle at every vertex is a rational multiple of $\pi$.
One can naturally associate the data of a surface with flat structure
$(S_{P},q_{P})$ so that the geodesics correspond to trajectories of
billiards in $P$.  We give a very brief discussion of this and refer the
reader to Zemlyakov and Katok \cite{ZK}, Kerckhoff, Masur, and Smillie
\cite{KMS}, and Masur and Tabachnikov \cite{MT} for more details.

To construct $S_{P}$, first consider the dihedral group $D_{2k}$ generated
by reflections in the lines through the origin in ${\mathbb R}^{2}$,
parallel to the sides of $P$.
Let
$${\cal P} = \coprod_{\gamma \in D_{2k}} \gamma P$$
be the disjoint union of the images of $P$ under the linear action by
$D_{2k}$.  We view the components $\gamma P$ as having a well-defined
embedding in ${\mathbb R}^{2}$, up to translation.  The group $D_{2k}$
acts on ${\cal P}$ in an obvious way, and we form the surface $S_{P}$
as the quotient of ${\cal P}$ obtained by identifying an edge $e$ of
${\cal P}$ with its image $\gamma e$, for $\gamma \in D_{2k}$, if $e$
and $\gamma e$ are parallel (with the same orientation).  $dz^{2}$
is defined on each polygon and pieces together to give a well-defined
quadratic differential $q_{P}$ on $S_{P}$.  In fact, $dz$ is invariant,
so that $q = \omega^{2}$ for a globally defined $1$--form $\omega$.
The polygon $P$ is said to be a {\em lattice polygon} if $(S_{P},q_{P})$
defines a Teichm\"uller curve in ${\cal M}(S_{P})$.

The right and isosceles lattice triangles have been classified by Kenyon
and Smillie in \cite{KS}.  They conjectured that there were exactly
three non-isosceles, acute lattice triangles and proved this for a large
number of examples.  The conjecture was proven by Puchta in \cite{P}.
We collect these facts together in the following.

\begin{theorem}[Kenyon--Smillie, Puchta] \label{ksptheorem}
The right lattice triangles are those with smallest angle $\frac{\pi}{k}$,
$k \in {\mathbb Z}_{\geq 4}$.
The acute isosceles lattice triangles are those with smallest angle
$\frac{\pi}{k}$, $k \in {\mathbb Z}_{\geq 3}$.
There are precisely three acute, non-isosceles lattice triangles, namely
those with angles
$$(1) \quad \left( \frac{\pi}{4},\frac{\pi}{3},\frac{5 \pi}{12} \right)
  \quad \quad
(2) \quad \left( \frac{2 \pi}{9},\frac{\pi}{3},\frac{4 \pi}{9} \right)
  \quad \quad
(3) \quad \left( \frac{\pi}{5},\frac{\pi}{3},\frac{7 \pi}{15} \right)$$
\end{theorem}

We will discuss the surfaces and quadratic differentials associated to
these lattice triangles in more detail in Section \ref{whichgroupsection}.

\section{Matrices and graphs} \label{pfsection}

\subsection{Non-negative matrices} \label{nonnegmatsection}

Let $M$ be a square, $n \times n$ matrix with real entries.  The {\em
spectral radius of $M$} is the maximum of the moduli of its eigenvalues,
and we denote this by $\mu(M)$.

We say that $M$ is {\em non-negative} (respectively, {\em positive})
if the entries of $M$ are non-negative (respectively, positive) and in
this case we write $M \geq 0$ (respectively, $M > 0$).  Say that $M \geq
0$ is {\em irreducible} if for every $1 \leq i,j \leq n$ there is some
power, $M^{k}$, so that $(M^{k})_{ij} > 0$ (see \cite{G}).  If $M,M'
\geq 0$, then write $M \leq M'$ if $M_{ij} \leq M_{ij}'$ for every $1
\leq i,j \leq n$ and write $M < M'$ if in addition this inequality is
strict for {\em some} $1 \leq i,j \leq n$.  We similarly define $\vec{V}
\geq 0$, $\vec{V} > 0$, $\vec{V} \leq \vec{V}'$, and $\vec{V} < \vec{V}'$
for vectors $\vec{V}$ and $\vec{V}'$ in ${\mathbb R}^{n}$.

The following theorem on irreducible matrices will be useful (see \cite{G}
for a proof).

\begin{theorem}[Perron--Frobenius] \label{pf1}
Suppose that $M \geq 0$ is irreducible.
Then $M$ has a unique (up to scaling) non-negative eigenvector $\vec{V}$.
This vector is positive with eigenvalue $\mu = \mu(M) > 0$.
Moreover, for any non-negative vector $\vec{U} \neq 0$, we have
$$\min_{1 \leq i \leq n} \left(\frac{(M\vec{U})_{i}}{\vec{U}_{i}}\right)
  \leq \mu
  \leq \max_{1 \leq i \leq n} \left(\frac{(M\vec{U})_{i}}{\vec{U}_{i}}\right)$$
with either inequality being an equality if and only if $\vec{U}$ is a
multiple of $\vec{V}$.
\end{theorem}

{\bf Remark}\qua
We define $\frac{(M\vec{U})_{i}}{\vec{U}_{i}} = + \infty$ whenever
$\vec{U}_{i} = 0$.

When $M$ is irreducible, we refer to the eigenvalue $\mu(M)$ (equal to
the spectral radius) as the {\em Perron--Frobenius eigenvalue} (briefly,
PF eigenvalue) of $M$ and an associated eigenvector, as in the theorem,
is called a {\em Perron--Frobenius eigenvector} (briefly, PF eigenvector)
for $M$.

\subsection{Graphs} \label{graphsubsection}

For more details on spectral radii of graphs, see the survey article of
Cvetkovi\'c and Rowlinson \cite{CR}.  The author is thankful to Curt
McMullen for pointing out this reference which greatly simplified the
exposition.

Given any finite graph ${\cal G}$, one can associate to ${\cal G}$
a matrix, ${\cal A}d({\cal G})$, called the {\em adjacency matrix}, as
follows.  Labeling the vertices of ${\cal G}$ by $x_{1},\ldots,x_{n}$, the
$(i,j)$--entry of ${\cal A}d({\cal G})$ is defined to be the number of edges
connecting $x_{i}$ to $x_{j}$.  The {\em spectral radius of ${\cal G}$} is
defined to be $\mu ({\cal G}) = \mu({\cal A}d({\cal G}))$.  Note that when
${\cal G}$ is connected, ${\cal A}d({\cal G})$ is irreducible.  Indeed,
$(({\cal A}d({\cal G}))^{k})_{ij}$ is the number of combinatorial paths
of length $k$ from the $i$th vertex to the $j$th.  The following
is an elementary consequence of Theorem \ref{pf1}.

\begin{theorem} \label{pfgraph}
If ${\cal G}_{0} \subset {\cal G}$ is a subgraph of a connected graph
${\cal G}$, then $\mu({\cal G}_{0}) \leq \mu({\cal G})$, with equality
if and only if ${\cal G}_{0} = {\cal G}$.
\end{theorem}

From this theorem one easily obtains a proof of the following (which is
a special case of the classical result of Smith \cite{Sm}).

\begin{theorem}[Smith]
\label{smithlist}
The set of connected bipartite graphs ${\cal G}$ with $\mu({\cal G}) <
2$ are precisely the recessive graphs, and those with $\mu({\cal G}) =
2$ are precisely the critical graphs.
\end{theorem}

\begin{proof}
An explicit calculation (see eg \cite{CR}) shows that the spectral
radius of every critical graph is $2$.
Any connected bipartite graph ${\cal G}$ contains or is contained in
one of the critical graphs.
To see this, we note that if ${\cal G}$ is not a tree, then it contains
a cycle (of even length since ${\cal G}$ is bipartite).
Hence ${\cal P}_{2c} \subset {\cal G}$ for some $c$.
If ${\cal G}$ is a tree, then one of the following holds
\begin{enumerate} 
\item ${\cal G}$ is homeomorphic to an interval (and thus contained in
  some ${\cal Q}_{c}$),
\item ${\cal G}$ contains a vertex with valence at least $4$ (and so
  contains ${\cal Q}_{5}$),
\item ${\cal G}$ has at least two vertices of valence at least $3$
  (and so contains some ${\cal Q}_{c}$), or
\item ${\cal G}$ has exactly one vertex of valence $3$ and all other
  vertices of valence at most $2$.
\end{enumerate}
In case (4), by inspection, ${\cal G}$ is either contained in one of
${\cal Q}_{c}$, ${\cal R}_{7}$, ${\cal R}_{8}$, or ${\cal R}_{9}$,
or else it contains ${\cal R}_{7}$, ${\cal R}_{8}$, or ${\cal R}_{9}$.

The only connected proper subgraphs of the critical graphs are the
recessive graphs, and so any other connected graph contains some critical
graph.  The theorem now follows from Theorem \ref{pfgraph}.
\end{proof}

There is also a classification of graphs, similar to Smith's, for graphs
having spectral radius in the interval $(2,\sqrt{2 + \sqrt{5}}]$ due
to Cvetkovi\'c, Doob, and Gutman \cite{CDG} and Brouwer and Neumaier
\cite{BN}.  From this, we easily obtain the following.

\begin{theorem}[Cvetkovi\'c, Doob, Gutman, Brouwer, Neumaier]
\label{cdgbnmin}
Given any bipartite graph ${\cal G}$ with $\mu({\cal G}) > 2$, we have
$$\mu({\cal G}) \geq \mu_{L} \approx 2.0065936$$
with equality if and only if ${\cal G} = {\cal E}h_{10}$.
\end{theorem}

Here $\mu_{L}$ is the square root of the unique largest root of
\begin{equation} \label{muleqn}
x^{5} - 9x^4 + 27x^{3} - 31x^{2} +12x -1
\end{equation}
This polynomial is the square root of the characteristic polynomial for
the matrix $({\cal A}d({\cal E}h_{10}))^{2}$.

\begin{figure}[ht!]
\begin{center}
\psfig{file=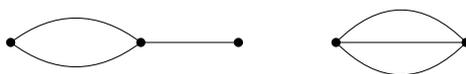,height=.7truein}
\caption{Subgraphs of a graph with multiple edges and $\mu > 2$}
\label{213}
\end{center}
\end{figure}

\begin{proof}
Appealing to the aforementioned classification (see \cite{CR}), one can
verify by explicit calculation that ${\cal E}h_{10}$ uniquely minimizes
spectral radius among graphs in the list, and that its spectral radius
is $\mu_{L}$.
The classification is for graphs without multiple edges, so we verify
directly that a graph ${\cal G}$ with multiple edges has $\mu({\cal G})
> \mu_{L}$.
Such a graph must contain one of the graphs shown in Figure \ref{213}.
These each have spectral radius at least $\sqrt{5} > \mu_{L}$, and so
the theorem follows from Theorem \ref{pfgraph}.
\end{proof}

\section{Affine actions for groups generated by two positive multi-twists} \label{singeuc}

\subsection{Constructing the flat structure} \label{curves2flatsection}

In this section, we recall the construction of Thurston \cite{T}.
Slight variations are also described in Long \cite{L}, Veech \cite{V1},
and a special case in the lecture notes \cite{FLP}.

Viewing $A \cup B$ as a graph on $S$, the components of $S \setminus A
\cup B$ are then the (interiors of) faces of this graph (actually, we
are viewing $S$ as a $2$--complex with $A \cup B$ as the $1$--skeleton).
Thus, each face is a disk (with at most one marked point) which we may
view as a $2k$--gon for some $k \in {\mathbb Z}$.  Since $A$ and $B$
are assumed to intersect minimally, any face containing no marked points
must have at least four edges.  Write $A = a_{1} \cup \cdots \cup a_{n}$
and $B = b_{1} \cup \cdots \cup b_{m}$.

\begin{figure}[ht!]\small
\centerline{}
\centerline{}
\begin{center}
\psfig{file=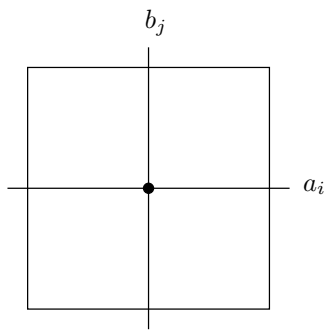,height=1.5truein}
\caption{The local picture in any rectangle}
\label{rectint}
\end{center}
  \setlength{\unitlength}{1in}
  \begin{picture}(0,0)(0,0)
    \put(3.4,1.25){$a_{i}$}
    \put(2.57,2.1){$b_{j}$}
  \end{picture}
\end{figure}

Let $\Gamma_{A,B}$ be the dual graph to $A \cup B$ embedded in $S$ so
that the vertex of $\Gamma_{A,B}$ dual to a face with a marked point is
that marked point.  $\Gamma_{A,B}$ defines a cell division of $S$, which
we also denote by $\Gamma_{A,B}$, each 2--cell of which is a rectangle.
Every rectangle contains a single arc of some $a_{i}$ and a single arc
of some $b_{j}$ intersecting in one point (see Figure \ref{rectint}).
Note that every vertex which is {\em not} a marked point of $S$ must
have valence at least $4$ by the previous paragraph.

One can now use $\Gamma_{A,B}$ to define a Euclidean cone metric on $S$
by declaring each rectangle to be a Euclidean rectangle.  The choice
of Euclidean rectangles is of course subject to the condition that
whenever two rectangles meet along an edge, the shared edge must have
the same length in each rectangle.  It follows that we obtain one real
parameter for each component of $A$ and of $B$, corresponding to the
length of the edges which that component meets.  This defines a flat
structure having orthonormal basis $\pm \{e_{1},e_{2}\}$ with $e_{1}$
parallel to the edges which $B$ transversely intersects, and $e_{2}$
parallel to the edges which $A$ intersects.

Since we want $\langle T_{A},T_{B} \rangle$ to act by affine
transformations with respect to this structure, we choose these rectangle
parameters as follows.
Define $N = N_{A,B}$ to be the $n \times m$ matrix whose $(i,j)$--entry
is $i(a_{i},b_{j})$.
The connectivity of $A \cup B$ guarantees that $NN^{t}$ is irreducible
(here $N^{t}$ is the matrix transpose of $N$).
Let $\vec{V}$ be a PF eigenvector for $\mu = \mu(NN^{t})$.
Notice that for the same reason, $N^{t}N$ is also irreducible, and
setting $\vec{V}' = \mu^{-\frac{1}{2}} N^{t} \vec{V} \geq 0$, we see that
$$N^{t}N \vec{V}' = N^{t} N \mu^{-\frac{1}{2}} N^{t} \vec{V} =
  \mu^{-\frac{1}{2}}N^{t} (N N^{t} \vec{V}) = \mu^{-\frac{1}{2}}N^{t}
  \mu \vec{V} = \mu \mu^{-\frac{1}{2}}N^{t} \vec{V} = \mu \vec{V}'$$
so that $\mu(N^{t}N) = \mu = \mu(NN^{t})$.  With this choice of $\vec{V}$
and $\vec{V}'$, note that we also have $\vec{V} = \mu^{-\frac{1}{2}}
N \vec{V}'$.  We write $\mu(A \cup B)$ to denote $\sqrt{\mu(NN^{t})}$
(the reason for the square root will soon become evident).

\begin{figure}[ht!]\small
\centerline{}
\centerline{}
\begin{center}
\psfig{file=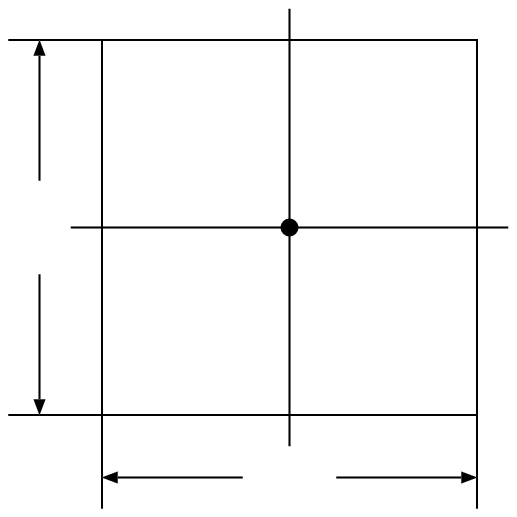,height=1.6truein}
\caption{The Euclidean rectangle}
\label{sides}
\end{center}
  \setlength{\unitlength}{1in}
  \begin{picture}(0,0)(0,0)
    \put(3.45,1.4){$a_{i}$}
    \put(2.67,2.2){$b_{j}$}
    \put(1.85,1.39){$\vec{V}_{i}$}
    \put(2.64,.57){$\vec{V}_{j}'$}
  \end{picture}
\end{figure}

We now make any rectangle of $\Gamma_{A,B}$ containing arcs of $a_{i}$
and $b_{j}$ into a Euclidean rectangle for which the sides transverse
to $a_{i}$ have length $\vec{V}_{i}$ and the sides transverse to $b_{j}$
have length $\vec{V}_{j}'$ (see Figure \ref{sides}).
For any component $a_{i}$ of $A$, the rectangles containing arcs of
$a_{i}$ fit together to give a Euclidean annulus (which is a neighborhood
of $a_{i}$).
The length of this annulus is $\vec{V}_{i}$, and to see what the girth
is, note that for each $j = 1,\ldots,m$ and for each intersection point
of $a_{i}$ with $b_{j}$ there is a rectangle of width $\vec{V}_{j}'$
in the annulus.
So, for each $j = 1,\ldots,m$, there is a contribution of $i(a_{i},b_{j})$
rectangles of width $\vec{V}_{j}'$.
Therefore, the girth is
$$\sum_{j=1}^{m} i(a_{i},b_{j}) \vec{V}_{j}' = (N \vec{V}')_{i} =
  \mu^{\frac{1}{2}} \vec{V}_{i} = \mu(A \cup B) \vec{V}_{i}$$
Similarly, the rectangles containing arcs of $b_{j}$ fit together to
give a Euclidean annulus of length $\vec{V}_{j}'$ and girth $\mu(A \cup
B) \vec{V}_{j}'$.

We now verify that $T_{A}$ and $T_{B}$ are represented by affine
transformations with respect to this structure.
The derivative of the affine map for $T_{A}$ (in terms of $\pm
\{e_{1},e_{2} \}$) is given by
$$\DAf(T_{A}) = \left( \begin{array}{cc}
  1 & \mu(A \cup B) \\
  0 & 1\\ \end{array} \right)$$
To see this, first construct the affine twist on each of the Euclidean
annuli described above and note that it has the desired derivative (see
(\ref{dehneqn}) and Section \ref{multitwistsection}).
Since each of the twists is the identity on the boundary of its
defining annulus, they all piece together to give a well-defined affine
homeomorphism with the correct derivative.
Similarly, the derivative of the affine representative of $T_{B}$ is
$$\DAf(T_{B}) = \left( \begin{array}{cc}
  1 & 0 \\
  -\mu(A \cup B) & 1 \\ \end{array} \right)$$
We note that, by construction, all vertices of $\Gamma_{A,B}$ are fixed
by every element of $\langle T_{A},T_{B} \rangle$.

\subsection{$\mu(A \cup B)$ vs $\mu({\cal G}(A \cup B))$} \label{muvmusect}

Given $A \cup B$ filling $S$, we have associated two positive numbers,
$\mu(A \cup B)$ and $\mu({\cal G}(A \cup B))$.
Not surprisingly, these are the same numbers.

\begin{proposition} \label{whatismu}
With $N = N_{A,B}$ as in the previous section, we have
$${\cal A}d({\cal G}(A \cup B)) =
  \left( \begin{array} {cc}
  0 & N \\
  N^{t} & 0 \\ \end{array} \right)$$
In particular, $\mu(A \cup B) = \sqrt{\mu(NN^{t})} = \mu({\cal G}(A \cup B))$.
\end{proposition}

\proof
Let us denote the vertices of ${\cal G}(A \cup B)$ (and curves of $A$ and
$B$) by both $a_{1},\ldots,a_{n},b_{1},\ldots,b_{m}$ and $x_{1},\ldots,x_{n+m}$
where $x_{i} = a_{i} \subset A$, $i=1,\ldots,n$, and $x_{j+n} = b_{j}
\subset B$, $j=1,\ldots,m$.
Then, the $(i,j)$--entry of ${\cal A}d({\cal G}(A \cup B))$ is the number
of edges from $x_{i}$ to $x_{j}$.
By definition of ${\cal G}(A \cup B)$, this is equal to $i(x_{i},x_{j})$.
Since the $(i,j)$--entry of $N$ (respectively $N^{t}$) is $i(a_{i},b_{j})$
(respectively $i(b_{i},a_{j})$), it is immediate that
$${\cal A}d({\cal G}(A \cup B)) =
  \left( \begin{array} {cc}
  0 & N \\
  N^{t} & 0 \\ \end{array} \right)$$
To see the second statement, note that $\mu({\cal A}d({\cal G}(A \cup
B)))^{2} = \mu(({\cal A}d({\cal G}(A \cup B)))^{2})$ and that
$$({\cal A}d({\cal G}(A \cup B)))^{2} = 
  \left( \begin{array} {cc}
  NN^{t} & 0 \\
  0 & N^{t}N \\ \end{array} \right)\eqno{\qed}$$

\section{Fuchsian groups} \label{fuchsian}

Here we note a few lemmas concerning the Fuchsian groups which occur
as the images of groups $\langle T_{A},T_{B} \rangle$ under $\DAf$.
The proofs are routine exercises in hyperbolic geometry, and we refer
to Beardon's text \cite {Be}, Ratcliffe's text \cite{Ra}, and the
Thurston's notes \cite{Tn} for more background on hyperbolic geometry
and Fuchsian groups.  The following two theorems are then easily derived
from these lemmas.

\begin{theorem} \label{main1}
$\langle T_{A},T_{B} \rangle \cong {\mathbb F}_{2}$ if and only if ${\cal
G}(A \cup B)$ contains a dominant component.
\end{theorem}

\begin{theorem} \label{main5}
For any surface $S$, any $A, B \in {\cal S}'(S)$, and any pseudo-Anosov
element
$$\phi \in \langle T_{A},T_{B} \rangle <  \Mod(S)$$
we have $\lambda(\phi) \geq \lambda_{L} \approx 1.1762808$.
Moreover, $\lambda(\phi) = \lambda_{L}$ precisely when $S$ has genus
$5$ (with at most one marked point), $\{ A,B \} = \{ A_{L}, B_{L} \}$
as in Figure \ref{lehmerscurves} (up to homeomorphism), and $\phi$
is conjugate to $(T_{A}T_{B})^{\pm 1}$.
\end{theorem}

{\bf Remark}\qua I would like to thank Joan Birman for pointing out that
the number obtained as the minimal dilation was Lehmer's number.

\subsection{Groups generated by two parabolics} \label{2parabolicsection}

For $\mu > 0$, set
$$\gamma_{1}(\mu) = \left( \begin{array}{cc}
  1 & \mu \\
  0 & 1 \\ \end{array} \right), \quad
  \gamma_{2}(\mu) = \left( \begin{array}{cc}
  1 & 0 \\
  -\mu & 1\\ \end{array} \right) \in \PSL_{2}{\mathbb R}$$
and
$$\Gamma_{\mu} = \, \langle \gamma_{1}(\mu), \gamma_{2}(\mu)  \rangle <
  \PSL_{2}{\mathbb R}$$
Note that $\DAf$ maps $\langle T_{A},T_{B} \rangle$ onto $\Gamma_{\mu(A
\cup B)}$.  We write ${\cal O}_{\mu} = {\mathbb H}^{2}/\Gamma_{\mu}$.

Recall that the convex core of a hyperbolic manifold $M$, denoted $C(M)$,
is the smallest convex sub-manifold for which the inclusion is a homotopy
equivalence.
The {\em signature} of a triangle orbifold, ${\cal O}$, (or, equivalently,
its associated triangle group) is a triple $(p,q,r)$ with $p,q,r \in
{\mathbb Z}_{\geq 2} \cup \{ \infty \}$.
${\cal O}$ is a sphere with cone points of order $p$, $q$, and $r$,
where a cone point of infinite order is a cusp.

\begin{lemma} \label{fuchtheorem}
If $\Gamma_{\mu}$ is discrete, then $\Gamma_{\mu} \cong {\mathbb F}_{2}$
if and only if $\mu \geq 2$.  Moreover,
\begin{enumerate}
\item for $\mu > 2$, ${\cal O}_{\mu}$ has infinite area and
  $C({\cal O}_{\mu})$ is a twice punctured disk,
\item for $\mu = 2$, ${\cal O}_{\mu} \cong \dot{S}_{0,3}$, and
\item for $\mu < 2$, ${\cal O}_{\mu}$ is a triangle orbifold with signature
  \begin{itemize}
  \item $(q,\infty,\infty)$ if $\gamma_{1}$ and $\gamma_{2}$ are not
    conjugate and
  \item $(2,q,\infty)$ if $\gamma_{1}$ and $\gamma_{2}$ are conjugate 
  \end{itemize}
where $q$ is the order of $\gamma_{1}\gamma_{2}$.
\end{enumerate}
In particular, $\Gamma_{\mu}$ has finite co-area if and only if $\mu \leq 2$.
\end{lemma}

\begin{proof}
If $\mu \geq 2$, one can construct a fundamental domain for the
action of $\Gamma_{\mu}$ on ${\mathbb H}^{2}$, as is shown in Figure
\ref{fundomain1} in the upper half-plane model.
Identifying the faces of this fundamental domain as indicated gives the
quotient ${\cal O}_{\mu}$.
Therefore ${\cal O}_{\mu} \cong \dot{S}_{0,3}$ for $\mu = 2$, and $C({\cal
O}_{\mu})$ is a twice-punctured disk (with punctures represented by cusps)
for $\mu > 2$.
In particular, $\Gamma_{\mu} \cong {\mathbb F}_{2}$ if $\mu \geq 2$.

\begin{figure}[ht!]\small
\begin{center}
\psfig{file=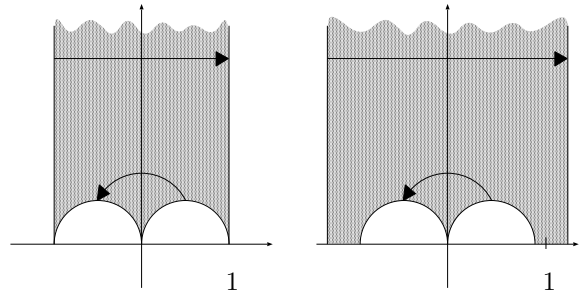,height=1.5truein}
\caption{Fundamental domains for $\Gamma_{\mu}$ with $\mu = 2$ (left)
and $\mu > 2$ (right)}
\label{fundomain1}
\end{center}
  \setlength{\unitlength}{1in}
  \begin{picture}(0,0)(0,0)
    \put(2.23,.6){$1$}
    \put(3.89,.6){$1$}
  \end{picture}
\end{figure}

When $\mu < 2$, we note that $\Tr(\gamma_{1}\gamma_{2}) = 2 - \mu^{2}
\in (-2,2)$, and so $\gamma_{1}\gamma_{2}$ is elliptic.
Because $\Gamma_{\mu}$ is discrete, $\gamma_{1}\gamma_{2}$ must have
finite order, and so $\Gamma_{\mu} \not \cong {\mathbb F}_{2}$.

Using the Dirichlet domain construction centered at the point $2i$,
one can check that the resulting fundamental domain must be contained
in the set, $P$, shown in Figure \ref{fundomain2}.
$$P = \left\{ z \in {\mathbb H}^{2} \, | \, d(z,2i) \leq d(z,\gamma(2i))
  \mbox{ for } \gamma = \gamma_{1}(\mu)^{\pm 1},
  \gamma_{2}(\mu)^{\pm 1} \right\}$$
$P$ is a polygon with two finite vertices at points
$z_{\pm} \in {\mathbb H}^{2}$ ($Re(z_{\pm}) = \frac{\pm \mu}{2}$
and $|z_{\pm}| = 1$) and two infinite vertices at $0$ and $\infty$.

\begin{figure}[ht!]\small
\begin{center}
\psfig{file=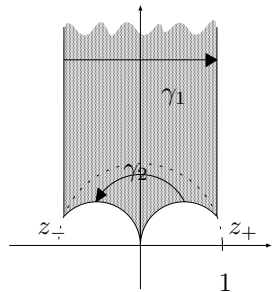, height=1.5truein}
\caption{The set $P$ containing the fundamental domain}
\label{fundomain2}
\end{center}
  \setlength{\unitlength}{1in}
  \begin{picture}(0,0)(0,0)
    \put(3,.6){$1$}
    \put(3.05,.9){$z_{+}$}
    \put(2.05,.9){$z_{-}$}
    \put(2.5,1.2){$\gamma_{2}$}
    \put(2.7,1.6){$\gamma_{1}$}
  \end{picture}
\end{figure}

The point $z_{+}$ is fixed by $\gamma_{1}\gamma_{2}$ and $z_{-}$ is fixed
by the conjugate $\gamma_{2}\gamma_{1}$.  Let $\theta$ denote the interior
angle at $z_{+}$, which is equal to the angle at $z_{-}$.  A calculation
shows that $2\cos\left(\frac{\theta}{2}\right) = \mu$ and that
$\gamma_{1}\gamma_{2}$ is a {\em clockwise} rotation about the point
$z_{+}$ through an angle $2 \theta$.

If $2 \theta = \frac{2 \pi}{q}$ for some integer $q$, then ${\cal O}$
is obtained from $P$ by identifying the two pairs of edges according
to $\gamma_{1}$ and $\gamma_{2}$.  In this case, ${\cal O}_{\mu}$ has
two cusps, one corresponding to each of $\langle \gamma_{1} \rangle$
and $\langle \gamma_{2} \rangle$ and one cone point of order $q =
\frac{2 \theta}{2 \pi}$ (which is the order of $\gamma_{1}\gamma_{2}$).
The maximal parabolic subgroups $\langle \gamma_{1} \rangle$ and
$\langle \gamma_{2} \rangle$ are not conjugate since they represent
different cusps.  In this case ${\cal O}_{\mu}$ is a triangle orbifold
with signature $(q,\infty,\infty)$, as required.

If $2 \theta = \frac{4 \pi}{q}$ for some odd integer $q$, then
$\gamma_{1}\gamma_{2}$ generates a cyclic subgroup of order $q$ in
$\Gamma_{\mu}$.  This subgroup also contains the clockwise rotation
$\rho = (\gamma_{1}\gamma_{2})^{\frac{q+1}{2}}$ about $z_{+}$ through
an angle $\theta = \frac{2 \pi}{q}$.  Consider the element $\delta
= \gamma_{1}^{-1}\rho$ which takes $0$ to $\infty$ and fixes $i$.
If we intersect $P$ with the complement of the unit disk in ${\mathbb
C}$, we obtain a fundamental domain for $\Gamma_{\mu}$ as shown in
Figure \ref{fundomain3}, with quotient given by the identifications
by $\gamma_{1}$ and $\delta$ as shown.  In this case ${\cal O}_{\mu}$
is a triangle orbifold with signature $(2,q,\infty)$ and $\gamma_{1}$
and $\gamma_{2}$ are conjugate by $\delta$, as required.

\begin{figure}[ht!]\small
\begin{center}
\psfig{file=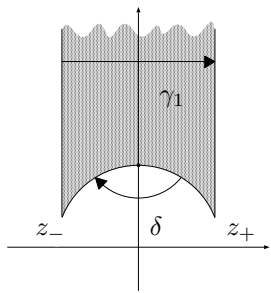, height=1.5truein}
\caption{The fundamental domain when $2 \theta = \frac{4 \pi}{q}$, $q$ odd}
\label{fundomain3}
\end{center}
  \setlength{\unitlength}{1in}
  \begin{picture}(0,0)(0,0)
    \put(3.05,.9){$z_{+}$}
    \put(2.05,.9){$z_{-}$}
    \put(2.7,1.6){$\gamma_{1}$}
    \put(2.65,.9){$\delta$}
  \end{picture}
\end{figure}

Finally, suppose $2 \theta \not \in \left\{ \frac{4\pi}{q} \,\, | \,\,
q \in {\mathbb Z}_{+} \right\}$.
In this case the cyclic subgroup generated by $\gamma_{1}\gamma_{2}$
contains a rotation about $z_{+}$ through an angle less than $\theta$.
In particular, the Dirichlet domain construction based at an appropriate
point on the imaginary axis gives a {\em compact} fundamental domain
contained in $P$.
This contradicts the fact that $\Gamma_{\mu}$ contains parabolics and
hence any fundamental domain is noncompact.
\end{proof}

\subsection{Freeness at last} \label{freenessproofsection}

{\bf Proof of Theorem \ref{main1}}\qua
Proposition \ref{reductionproposition} shows that it suffices to
assume that $A \cup B$ fills $S$.  According to Theorem \ref{smithlist},
$\mu({\cal G}(A \cup B)) \geq 2$ if and only if ${\cal G}(A \cup B)$ is
dominant.  Now, by Proposition \ref{whatismu}, $\mu(A \cup B) = \mu({\cal
G}(A \cup B))$, so it suffices to show that $\langle T_{A},T_{B} \rangle$
is free if and only if $\mu(A \cup B) \geq 2$.  This latter equivalence
is precisely what Thurston had suggested in \cite{T}.

Suppose $\mu(A \cup B) \geq 2$.
By Lemma \ref{fuchtheorem}, $\DAf$ is a surjection to ${\mathbb F}_{2}$.
Since $\langle T_{A},T_{B} \rangle$ is generated by two elements, we
have a surjection
$$\chi\co {\mathbb F}_{2} \rightarrow \langle T_{A},T_{B} \rangle$$
$\DAf \circ \chi$ is therefore a surjection from ${\mathbb F}_{2}$
onto itself.
Free groups are Hopfian (see \cite{MKS}), hence this is an isomorphism.
Therefore, $\chi$ is also an isomorphism, proving that $\langle
T_{A},T_{B} \rangle \cong {\mathbb F}_{2}$.

If $\mu(A \cup B) < 2$, Lemma \ref{fuchtheorem} implies $\DAf(T_{A}T_{B})$
has finite order, so Theorem \ref{reptheorem} says $T_{A}T_{B}$ also
has finite order, hence $\langle T_{A},T_{B} \rangle \not \cong
{\mathbb F}_{2}$.
\endproof

The next proposition is included for it own interest.  It makes precise
the statement that most elements in $\langle T_{A},T_{B} \rangle$
are pseudo-Anosov.

\begin{proposition} \label{elementtypes}
Suppose $A \cup B$ fills $S$ and ${\cal G}(A \cup B)$ is dominant.
Then every element of $\langle T_{A},T_{B} \rangle$ is pseudo-Anosov
except conjugates of powers of $T_{A}$ and $T_{B}$, and also $T_{A}T_{B}$
when ${\cal G}(A \cup B)$ is critical.
\end{proposition}

\begin{proof}
The hyperbolic elements of $\Gamma_{\mu} = \pi_{1}({\mathbb
H}^{2}/\Gamma_{\mu})$ are precisely those elements corresponding to
loops that are freely homotopic to closed geodesics.  All loops that are
not homotopic to cusps have such a representative, and so the corollary
follows from Theorem \ref{reptheorem}, Lemma \ref{fuchtheorem}, and the
fact that $\DAf$ is an isomorphism in this case.
\end{proof}

{\bf Remark}\qua
When $\Gamma_{\mu}$ is discrete and $\mu < 2$, it also has precisely $3$
conjugacy classes of cyclic subgroups which make up all non-hyperbolic
elements.  However, $\DAf$ is not necessarily an isomorphism in this
case, so we only know that all non-pseudo-Anosov elements map by $\DAf$
to one of $3$ cyclic subgroups, up to conjugacy.

\subsection{Dilatation bounds and Lehmer's number}
\label{translationlengthboundssection}

The connection between translation lengths and dilatations is provided
by Theorem \ref{reptheorem}.  This is the basis for the proof of Theorem
\ref{main5}, and so to apply this we will need a few elementary facts
concerning translation lengths in the Fuchsian groups $\Gamma_{\mu}$.
As Lemma \ref{fuchtheorem} shows, for $\mu \leq 2$, ${\mathbb H}^{2} /
\Gamma_{\mu}$ is a triangle orbifold.  Furthermore, triangle orbifolds
can have no closed embedded geodesics, so the following is an immediate
consequence of Theorem 11.6.8 of Beardon \cite{Be}.

\begin{proposition}[Beardon]
\label{shortesttrianglegeodesics}
For any $\mu \leq 2$, the smallest translation length of a hyperbolic
element of $\Gamma_{\mu}$ is bounded below by
$$2\sinh^{-1} \left( \sqrt{\cos\left(\frac{3 \pi}{7}\right)} \right)$$
\end{proposition}

For the remaining cases, we have the following.

\begin{lemma} \label{othershortgeodesics}
When $\mu > 2$, the smallest translation length of a hyperbolic element of
$\Gamma_{\mu}$ is realized (uniquely up to conjugacy) by $(\gamma_{1}(\mu)
\gamma_{2}(\mu))^{\pm 1}$, and is given by $2 \log (\lambda_{\mu})$,
where $\lambda_{\mu}$ is the larger root of
\begin{equation} \label{lambdaeqn}
x^{2} + x(2 - \mu^{2}) + 1
\end{equation}
\end{lemma}

{\bf Remark}\qua
The larger root $\lambda_{\mu}$ of (\ref{lambdaeqn}) defines an increasing
function of $\mu > 2$.\\

\begin{proof}
Since $C({\mathbb H}^{2}/ \Gamma_{\mu})$ is a twice-punctured disk,
any geodesic $\gamma$ determines a conjugacy class represented by an
element which we also call $\gamma \in \Gamma_{\mu}$.  The translation
length of $\gamma$ is the length of the geodesic with the same name, and
is given by $2 \log(\lambda)$ where $\lambda$ is the spectral radius of
a matrix representative of $\gamma \in \Gamma_{\mu}$.  The only simple
closed geodesic is the boundary of the convex core.  Moreover, for
any other closed geodesic, one can cut and paste a collection of arcs
of this geodesic to obtain a curve homotopic to this boundary curve.
It follows that the boundary geodesic is the unique shortest geodesic.
This is represented by $\gamma_{1}(\mu)\gamma_{2}(\mu)$.

The natural representation of the projective class of
$\gamma_{1}(\mu)\gamma_{2}(\mu)$ by a matrix (given the matrices
for $\gamma_{1}(\mu)$ and $\gamma_{2}(\mu)$ we have chosen) has
$Tr(\gamma_{1}(\mu)\gamma_{2}(\mu)) = 2 - \mu^{2} < 0$, so we see that
$-\lambda$ satisfies the characteristic equation.  Therefore, $\lambda$
is the larger root of (\ref{lambdaeqn}).
\end{proof}

\begin{corollary} \label{dilboundcor}
Let $\phi \in \langle T_{A},T_{B} \rangle$ be any pseudo-Anosov
automorphism.  If $\mu = \mu(A \cup B) > 2$, then
$$\lambda(\phi) \geq \lambda_{\mu}$$
where $\lambda_{\mu}$ is a root of (\ref{lambdaeqn}).
Equality holds if and only if $\phi = (T_{A}T_{B})^{\pm 1}$ up to
conjugacy.  If $\mu \leq 2$, then
$$\lambda(\phi) > 1.47$$
\end{corollary}

{\bf Proof}\qua
Suppose $\mu > 2$.
By Theorem \ref{reptheorem} $\lambda(\phi) =
\exp(\frac{1}{2}L(\DAf(\phi)))$, hence the smallest dilatation occurs
precisely when $\DAf(\phi)$ has smallest translation length (dilation
is an increasing function of translation length).
By Lemma \ref{othershortgeodesics}, this is precisely when $\DAf(\phi)$
is conjugate to $\DAf(T_{A}T_{B})^{\pm 1}$.
As the proof of Theorem \ref{main1} shows, $\DAf$ is an isomorphism, and
so this happens if and only if $\phi$ is conjugate to $(T_{A}T_{B})^{\pm 1}$.
In this case, we have $\lambda(\phi) = \lambda(T_{A}T_{B}) =
\exp\left(\frac{1}{2} 2 \log(\lambda_{\mu})\right) = \lambda_{\mu}$.

When $\mu \leq 2$, by similar reasoning (appealing now to
Proposition \ref{shortesttrianglegeodesics} rather than Lemma
\ref{othershortgeodesics}) we obtain
$$\lambda(\phi) \geq \exp \left( \frac{1}{2} 2\sinh^{-1} \left(
  \sqrt{\cos\left(\frac{3 \pi}{7}\right)} \right) \right) > 1.47 \eqno{\qed}$$
{\bf Proof of Theorem \ref{main5}}\qua
By Corollary \ref{dilboundcor}, we need only consider
$\mu = \mu(A \cup B)> 2$, and it suffices to show that
$\lambda_{\mu} \geq \lambda_{L}$
with equality if and only if $\{A,B\} = \{A_{L},B_{L}\}$.
The remark following Lemma \ref{othershortgeodesics} tells us that to
minimize $\lambda_{\mu}$, we need only minimize $\mu(A \cup B)$.
Theorem \ref{cdgbnmin} says that $\mu$ is minimized uniquely by $\mu_{L}$
when ${\cal G}(A \cup B) = {\cal E}h_{10}$.
Here $\mu_{L}^{2}$ is the largest root of the polynomial (\ref{muleqn}).
By Corollary \ref{dilboundcor}, $\lambda_{\mu_{L}}$ is a root of
(\ref{lambdaeqn}) with $\mu = \mu_{L}$.
Thus, $(\lambda_{\mu_{L}},\mu_{L}) = (x,y)$ satisfies
$$\left\{ \begin{array}{l}
  x^{2} + x(2 - y^{2}) + 1 = 0 \\
  y^{5} - 9y^4 + 27y^{3} - 31y^{2} +12y -1 = 0\\ \end{array}\right.$$
Eliminating $y$ from this pair we find that $\lambda_{\mu_{L}}$ is a
root of (\ref{lehmerseqn}), and so $\lambda_{\mu_{L}} = \lambda_{L}$.

The only configuration which minimizes $\mu$ is ${\cal E}h_{10}$.
Since ${\cal E}h_{10}$ is a tree with one vertex having valence at most
three, Proposition \ref{recessiveuniqueness} completes the proof of the
theorem.
\endproof

\section{Teichm\"uller curves, triangle groups, and billiards}
\label{whichgroupsection}

The first theorem of this section is the following:

\begin{theorem} \label{main2}
The Teichm\"uller curves for which the associated stabilizers contain
a group generated by two positive multi-twists with finite index are
precisely those defined by $A \cup B$ filling $S$, where ${\cal G}(A
\cup B)$ is critical or recessive.
\end{theorem}

\begin{proof}
We must show that $\langle T_{A},T_{B} \rangle$ has finite co-area if
and only if ${\cal G}(A \cup B)$ is recessive or critical.
Proposition \ref{mcproposition} implies that the former happens if and only if $\Gamma_{\mu(A \cup B)} = \langle \DAf(T_{A}) ,\DAf(T_{B}) \rangle$ has finite co-area.
By Lemma \ref{fuchtheorem}, $\Gamma_{\mu(A \cup B)}$ has finite co-area
if and only if $\mu(A \cup B) \leq 2$.
Proposition \ref{whatismu} implies $\mu(A \cup B) = \mu({\cal G}(A
\cup B))$, and Theorem \ref{smithlist} implies that $\mu({\cal G}(A
\cup B)) \leq 2$ if and only if ${\cal G}(A \cup B)$ is recessive or
critical.
\end{proof}

The Teichm\"uller curves obtained from this theorem are most interesting
when ${\cal G}(A \cup B)$ are recessive (the critical configurations
all give Teichm\"uller curves covered by the thrice-punctured sphere).
To better understand these curves, we describe another construction
for surfaces and flat structures studied by Veech \cite{V1} and Earle
and Gardiner \cite{EG}.  Embedded in this construction is the billiard
construction for all but three of the lattice triangles from Theorem
\ref{ksptheorem}.  To complete the billiard picture we describe the
remaining three exceptional triangles and verify the following.

\begin{theorem} \label{main3}
The Teichm\"uller curves determined by the right and acute lattice
triangles have associated stabilizers containing a finite index subgroup
of the form $\langle T_{A},T_{B} \rangle$ with ${\cal G}(A \cup B)$
recessive.
\end{theorem}

The Teichm\"uller curves determined by these lattice triangles do
not account for all Teichm\"uller curves determined by recessive
configuration.  However, the constructions described below are general
enough to take care of all of these.  From this we obtain a complete
description of the non-free groups generated by two positive multi-twists.

\begin{theorem} \label{main4}
If ${\cal G}(A \cup B)$ is recessive, then $\DAf$ maps $\langle
T_{A},T_{B} \rangle$ onto a Fuchsian triangle group with finite central
kernel of order at most $2$.
The signature of the triangle group is described by the following table.
\def\strutt{\vrule width 0pt height 15pt depth 5pt}
$$\begin{array} {cccc}
\strutt  \mbox{configuration graph} & \mbox{signature} & \mbox{configuration graph}
    & \mbox{signature} \\
  \hline
  \strutt{\cal D}_{c}, \, c \geq 4 & (c-1,\infty,\infty) & {\cal E}_{6} &
    (6,\infty,\infty)\\
  \strutt{\cal A}_{2c+1}, \, c \geq 1 & (c+1,\infty,\infty) & {\cal E}_{7} &
    (9,\infty,\infty)\\
  \strutt{\cal A}_{2c}, \, c \geq 1 & (2,2c+1,\infty) & {\cal E}_{8} &
    (15,\infty,\infty)\\
  \end{array}$$
\end{theorem}

We will explain the proof of Theorems \ref{main3} and \ref{main4} in
Section \ref{outlineproofsection}.
The strategy in all cases is the same, but requires verification on a
case-by-case basis.
This is done for the ${\cal D}_{c}$ graphs in Section \ref{1hangersection},
the ${\cal A}_{c}$ graphs in Section \ref{oddstraightchainsection} for $c$
odd and Section \ref{evenstraightchainsection} for $c$ even, and finally
the three graphs ${\cal E}_{6}$, ${\cal E}_{7}$, and ${\cal E}_{8}$
in Section \ref{threeexceptionsection}.

\subsection{Proof outlined} \label{outlineproofsection}

In sections \ref{1hangersection} through \ref{threeexceptionsection} we
will describe a flat structure and positive multi-twists in $1$--manifolds
$A$ and $B$ which act as affine automorphisms with respect to this
structure.  This will generally be a different structure than the one
we constructed in Section \ref{curves2flatsection} for $A$ and $B$.
However, it is affine equivalent to that one, ie they define the
same Teichm\"uller disk.  This follows from the uniqueness in Theorem
\ref{bersuniquetheorem} since both Teichm\"uller disks are stabilized
by any pseudo-Anosov automorphism in $\langle T_{A},T_{B} \rangle$.

For each of the flat structures under consideration, and each of the
pairs of $1$--manifolds $A$ and $B$, we will see that ${\cal G}(A \cup B)$
is recessive.  In particular, to prove Theorem \ref{main3} it suffices
to recognize those flat structures coming from the billiard construction
for the right and acute lattice triangles among those which we describe.
This is verified in Section \ref{threeexceptionsection}.

\medskip{\bf Proof of Theorem \ref{main4}}\qua
$\DAf$ maps $\langle T_{A},T_{B} \rangle$ onto a Fuchsian triangle
group by Lemma \ref{fuchtheorem} and the fact that $\mu(A \cup B)
< 2$.  To see that the kernel is central, we note that by definition,
any element of the kernel has derivative $\pm I$.  In particular, this
must leave both the $A$--annuli and the $B$--annuli invariant, and hence
also each of $A$ and $B$ are invariant.  Any automorphism which leaves
a $1$--manifold invariant must commute with the associated multi-twist.
In particular, we see that every element of the kernel commutes with
the generators $T_{A}$ and $T_{B}$, and so is central.

Since any element of the kernel leaves each of $A$ and $B$ invariant,
it induces an automorphism of the graph $A \cup B$ which leaves the $A$
edges and $B$ edges invariant.  This in turn induces an automorphism of
the graph ${\cal G}(A \cup B)$ preserving the bicoloring.
Said differently, we obtain a homomorphism
$$\delta\co \ker(\DAf) \rightarrow \Aut_{bc}({\cal G}(A \cup B))$$
Here $\Aut_{bc}({\cal G}(A \cup B))$ is the automorphism group of the
graph preserving the bicoloring (which has index at most $2$ in the full
automorphism group).

\medskip
{\bf Claim}\qua\sl
$\ker(\delta)$ has order at most $2$.\rm\medskip

{\bf Proof of Claim}\qua
We note that any $\phi \in \ker(\delta)$ leaves each component of $A$
and of $B$ invariant.  Relabeling $A$ and $B$ and renumbering the
components if necessary, we may assume that $a_{1} \subset A$ is a
component corresponding to a $1$--valent vertex of ${\cal G}(A \cup B)$.
This $a_{1}$ has only a single point of intersection with $(A \setminus
\{ a_{1} \}) \cup B$, and hence $a_{1}$ is the closure of an edge,
$e_{1}$, of the graph $A \cup B$.  Since $a_{1}$ is invariant by $\phi$,
so is $e_{1}$.  Thus, $\ker(\delta)$ consists of isotopy classes of
(orientation preserving) homeomorphisms of $S$ leaving the $1$--skeleton
($A \cup B$) of a cell structure invariant and fixing the edge $e_{1}$
(not necessarily pointwise).  Such a group has order at most two.\endproof

Now, when ${\cal G}(A \cup B)$ is of type ${\cal A}_{2c}$,
$\Aut_{bc}({\cal G}(A \cup B))$ is trivial, hence $\ker(\DAf)$ has order
at most $2$ by the claim.  When ${\cal G}(A \cup B)$ is of type ${\cal
A}_{2c+1}$, $c \geq 1$, $\Aut_{bc}({\cal G}(A \cup B))$ has order two.
However, in this case we claim that $\ker(\delta)$ is trivial.  This is
because the two possible elements in this group are the identity and a
hyperelliptic involution.  The latter is in the full stabilizer, but not
in $\langle T_{A},T_{B} \rangle$ because it does not fix the two vertices
of $\Gamma_{A,B}$.  This implies $\ker(\DAf)$ has order at most two in
this case also.

When ${\cal G}(A \cup B)$ is of type ${\cal D}_{c}$, $\ker(\delta)$ is
again trivial.  This is because on the curve corresponding to the valence
three vertex, any $\phi \in \ker(\delta)$ must fix the three points of
intersection with the curves corresponding to the three adjacent vertices.
It follows that $\phi$ is the identity on that curve, hence on all of $S$.
When $c \geq 5$, there is only one non-trivial automorphism of ${\cal
D}_{c}$, and so $\ker(\DAf)$ has order at most two, and we are done in
this case.  When $c = 4$, we again use the fact that all vertices of
$\Gamma_{A,B}$ are fixed to see that the $\ker(\DAf)$ is trivial.

For the three exceptional cases we note that $\Aut({\cal G}(A \cup B))$
is trivial when ${\cal G}(A \cup B) = {\cal E}_{7}$ or ${\cal E}_{8}$,
so $\ker(\DAf)$ has order at most two by the claim.  In the one remaining
case that ${\cal G}(A \cup B) = {\cal E}_{6}$, we note that $\Aut({\cal
G}(A \cup B))$ has order two.  However, the non-trivial element is induced
by an automorphism of the surface which does not fix the vertices of
$\Gamma_{A,B}$, hence is not in $\langle T_{A},T_{B} \rangle$.

All that remains is to verify that $\langle \DAf(T_{A}),\DAf(T_{B})
\rangle$ has the required signature.  We will check below that
$\DAf(T_{A}T_{B})$ has order given by the larger of the two finite
numbers listed in the signature.  For the cases of ${\cal G}(A \cup B)$
of type ${\cal D}_{c}$, ${\cal A}_{2c+1}$, ${\cal E}_{6}$, ${\cal E}_{7}$,
and ${\cal E}_{8}$ this will prove that the signature is as listed by
showing that there must be two cusps in these cases.  For then, the
signature is $(q,\infty,\infty)$, where $q$ is the order of the product
of the two parabolic generators by Lemma \ref{fuchtheorem}.

Suppose that in the cases listed there were only one cusp.  By Lemma
\ref{fuchtheorem} there is an element $\DAf(\phi)$ conjugating
$\DAf(T_{A})$ to $\DAf(T_{B})$.  Up to an element of the kernel, $\phi$
would conjugate $T_{A}$ to $T_{B}$.  Because the kernel is central with
order at most two, we obtain a conjugation of $T_{A}^{2}$ to $T_{B}^{2}$.
This cannot happen since this would imply a homeomorphism taking a union
of two copies of $A$ to a union of two copies of $B$ which is not possible
for the given configurations.

Finally, for the remaining cases ${\cal G}(A \cup B) = {\cal A}_{2c}$, we
find an element in $\langle \DAf(T_{A}),\DAf(T_{B}) \rangle$ conjugating
$\DAf(T_{A})$ to $\DAf(T_{B})$.  Lemma \ref{fuchtheorem} completes the
proof, modulo finding this conjugating element and verifying the orders
of $T_{A}T_{B}$.  This is carried out in the next four sections.

\subsection{The ${\cal D}_{c}$ configurations} \label{1hangersection}

The following is described in more detail by Earle and Gardiner in \cite{EG}.

Consider a regular $2k$--gon, $\Delta_{2k}$, with $k \in {\mathbb Z}_{\geq
2}$, embedded in the plane with two vertical edges.  Identifying opposite
edges by Euclidean translations we obtain a surface $S$ of genus $\left\lfloor
\frac{k}{2} \right\rfloor$.  Because the gluings are by isometry, we obtain
a Euclidean cone metric on $S$, and the foliation by horizontal lines
provides a holomorphic quadratic differential $q$ (this restricts to
$dz^{2}$ on $\Delta_{2k}$).  Let $\alpha_{2k} = \frac{\pi}{k} $ and
$\beta_{2k} = \frac{\alpha_{2k}}{2}$.

Note first that the counter-clockwise rotation about the center of
$\Delta_{2k}$ through an angle $\alpha_{2k}$ defines an isometry of $S$
of order $2k$.  We denote this by $\rho_{2k}$.

\begin{figure}[ht!]
\cl{
\psfig{file=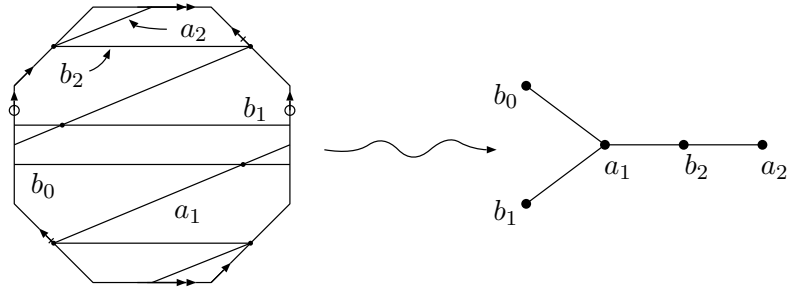,height=1.5truein}}
\caption{$A \cup B$ in $S$ and ${\cal G}(A \cup B)$, when $k = 4$}
\label{delta2kpoly}\vglue 4pt
  \setlength{\unitlength}{1in}
  \begin{picture}(0,0)(0,0)
    \put(1.5,.92){$a_{1}$}
    \put(1.53,1.87){$a_{2}$}
    \put(.75,1.05){$b_{0}$}
    \put(1.86,1.45){$b_{1}$}
    \put(.9,1.63){$b_{2}$}
    \put(3.75,1.15){$a_{1}$}
    \put(4.57,1.15){$a_{2}$}
    \put(3.17,1.52){$b_{0}$}
    \put(3.17,.91){$b_{1}$}
    \put(4.17,1.15){$b_{2}$}
  \end{picture}
\end{figure}

We also see that the horizontal foliation of $q$ has all closed leaves,
decomposing $S$ into $\left\lceil \frac{k}{2} \right\rceil$ annuli.
Let $B$ be the essential $1$--manifold which is the union of the cores
of the annuli, taking two parallel copies of the core of the
annulus meeting the two vertical sides of $\Delta_{k}$ (see Figure
\ref{delta2kpoly}).  $T_{B}$ acts by an affine transformation leaving
this foliation invariant, having derivative
$$\DAf(T_{B}) = \left( \begin{array}{cc}
  1 & 2 \cot (\beta_{2k}) \\
  0 & 1\\ \end{array} \right).$$

Rotate the horizontal foliation by an angle $\beta_{2k}$ (ie multiply
$q$ by $e^{i \alpha_{2k}}$).  This rotated foliation also has all closed
leaves, decomposing $S$ into $\left\lfloor \frac{k}{2} \right\rfloor$
annuli in another way.  Let $A$ be the union of the cores of these
annuli.  $T_{A}$ also acts by an affine transformation, with derivative
$\DAf(T_{A})$ given by
$$\left( \begin{array}{rr}
\cos(\beta_{2k}) & -\sin(\beta_{2k})\\
\sin(\beta_{2k}) & \cos(\beta_{2k})\\ \end{array} \right)
\left( \begin{array}{cc}
1 & 2 \cot(\beta_{2k}) \\
0 & 1 \\ \end{array} \right)
\left( \begin{array}{rr}
\cos(\beta_{2k}) & \sin(\beta_{2k})\\
-\sin(\beta_{2k}) & \cos(\beta_{2k})\\ \end{array} \right).$$
One can now verify that $T_{A}T_{B} = \rho_{2k}^{k+1}$ (eg consider
the action on the line segment from the center of the polygon to the
vertex at the top of the left vertical edge and on the horizontal line
segment from the center to the midpoint of the right vertical edge).
We also note that $\DAf(\rho_{2k}) \in \PSL_{2}{\mathbb R}$ has order $k$
(which is half its order in $\SL_{2}{\mathbb R}$).

As is indicated in Figure \ref{delta2kpoly}, ${\cal G}(A \cup B)$ is
the graph ${\cal D}_{k+1}$, when $k \geq 3$ (for $k = 2$, we get ${\cal
A}_{3}$).  This proves Theorem \ref{main4} for the graphs ${\cal D}_{c}$.

\subsection{The ${\cal A}_{c}$ configurations I: $c$ odd}
\label{oddstraightchainsection}

The examples below where studied by Veech in \cite{V1} and \cite{V2}.
However, we follow the discussion of Earle and Gardiner in \cite{EG}.

Let $S$, $A$, and $B$ be as in the previous section and assume that
$b_{0}$ and $b_{1}$ were the parallel components of $B$.  Write $A =
a_{1} \cup \cdots \cup a_{n}$ and $B = b_{0} \cup \cdots \cup b_{m}$.
Note that $n + m = k$.  Let $B' = B \setminus b_{0}$ and note that we
may replace $T_{B}$ by the following isotopic homeomorphism:
$$T_{B} \simeq \widehat{T}_{B'} = T_{b_{1}}^{2}T_{b_{2}} \cdots T_{b_{m}}$$
Now construct a $2$--fold cover $\pi\co\widetilde{S} \rightarrow S$ (which
is a branched cover when $k$ is odd) for which all components of $A$
and of $B' \setminus  b_{1} $ lift to loops, but the preimage of $b_{1}$
is a connected double cover of $b_{1}$.
Writing $\widetilde{A} = \pi^{-1}(A)$ and $\widetilde{B} = \pi^{-1}(B')$,
one can check that $T_{\widetilde{A}}$ and $T_{\widetilde{B}}$ cover
$T_{A}$ and $\widehat{T}_{B'}$, respectively.
Moreover, these act as affine transformations with respect to $\pi^{*}(q)$
with derivatives
$$\DAf(T_{\widetilde{A}}) = \DAf(T_{A}) \quad \mbox{ and } \quad \DAf(T_{\widetilde{B}}) = \DAf(\widehat{T}_{B'}) = \DAf(T_{B}).$$
So, $\DAf(T_{\widetilde{A}}T_{\widetilde{B}})$ has order $k$.

${\cal G}(\widetilde{A} \cup \widetilde{B})$ is of type ${\cal A}_{c}$
since each curve of $A$ and $B$ intersects at most two other curves.
$\widetilde{B}$ has $2m-1$ components, $\widetilde{A}$ has $2n$
components, so ${\cal G}(\widetilde{A} \cup \widetilde{B}) = {\cal A}_{2k
- 1}$.  This proves Theorem \ref{main4} for the graphs ${\cal A}_{c}$
with $c$ odd.

\medskip
{\bf Remark}\qua
In Veech's description of these examples, he explicitly constructed
the surface $\widetilde{S}$ from two regular $2k$--gons in the plane,
identified along an edge.  $\widetilde{S}$ is obtained by identifying
opposite sides of the resulting non-convex polygon.

\subsection{The ${\cal A}_{c}$ configurations II: $c$ even}
\label{evenstraightchainsection}

The following construction is due to Veech \cite{V1}, \cite{V2}.

For $k \in {\mathbb Z}_{\geq 1}$ we consider two regular $(2k+1)$--gons,
$\Delta_{2k+1}^{0}$ and $\Delta_{2k+1}^{1}$, in the plane sharing a
horizontal edge, and denote the non-convex polygon which is their union
by $\Delta_{2k+1}$.  Identifying opposite sides of $\Delta_{2k+1}$
we obtain a genus--$k$ surface, which we denote by $S$.

In the same fashion as above, we obtain a flat structure $q$ on $S$,
which restricts to $dz^{2}$ on $\Delta_{2k+1}$.  Let $\alpha_{2k+1} =
\frac{2\pi}{2k+1}$ and $\beta_{2k+1} = \frac{\alpha_{2k+1}}{2}$.

The counter-clockwise rotations through an angle $\alpha_{2k+1}$ about
the centers of $\Delta_{2k+1}^{0}$ and $\Delta_{2k+1}^{1}$ define
an isometry $\rho_{2k+1}$ of $S$ of order $2k+1$.  There is also an
involution $\sigma_{2k+1}$ obtained by rotating $\Delta_{2k+1}$ about the
center of the edge shared by $\Delta_{2k+1}^{0}$ and $\Delta_{2k+1}^{1}$.
Note that $\sigma_{2k+1}$ is in the kernel of $\DAf$.

The horizontal foliation has all closed leaves, and so decomposes $S$
into $k$ annuli.  Let $B$ be the union of the cores of these annuli.
Then $T_{B}$ acts on $S$ by affine transformations with derivative
$$\DAf(T_{B}) = \left( \begin{array}{cc}
1 & 2 \cot(\beta_{2k+1}) \\
0 & 1 \\ \end{array} \right)$$
Next, we let $A = \rho_{2k+1}^{k+1}(B)$.
Equivalently, $A$ is obtained as follows.
Rotate the horizontal foliation of $q$ through an angle $(k+1) \alpha_{2k+1}$.
This has the same effect as rotating through an angle $\beta_{2k+1} = (k+1)\alpha_{2k+1} - \pi$ (and hence multiplying $q$ by $e^{i \alpha_{2k+1}}$).
This foliation has all closed leaves and decomposes $S$ into annuli, the union of the cores of which are precisely $A$.

Now one can check that
$$T_{A}T_{B} = \rho_{2k+1} \sigma_{2k+1}$$
(eg one can verify that this holds on appropriately chosen segments).
So that we see
$$\DAf(T_{A}) = \DAf(T_{A}T_{B})^{k+1} \DAf(T_{B}) \DAf(T_{A}T_{B})^{-(k+1)}.$$
Thus, $\DAf(T_{A})$ and $\DAf(T_{B})$ are conjugate in $\DAf(\langle
T_{A},T_{B} \rangle)$ and $\DAf(T_{A}T_{B})$ has order $2k + 1$.  One can
check that ${\cal G}(A \cup B) = {\cal A}_{2k}$, thus proving Theorem
\ref{main4} for this class of graphs.

\subsection{Billiards and ${\cal E}_{6}$, ${\cal E}_{7}$, and ${\cal E}_{8}$}
\label{threeexceptionsection}

The constructions of the previous three sections provide a description
of the surfaces and quadratic differentials coming from the billiard
construction (see Section \ref{billiardsection}) for the right and acute
isosceles lattice triangles given in Theorem \ref{ksptheorem}, as we
shall now explain.

\begin{figure}[ht!]\small
\cl{
\psfig{file=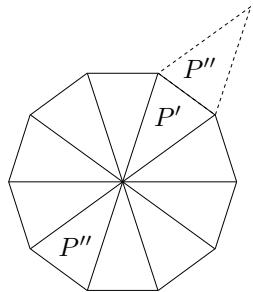,height=1.5truein}}
\caption{Reflecting $P'$ in the side opposite the apex gives a translate of $P''$}
\label{billiardpolyreflect}\vglue 4pt
  \setlength{\unitlength}{1in}
  \begin{picture}(0,0)(0,0)
    \put(2.74,1.4){$P'$}
    \put(2.89,1.65){$P''$}
    \put(2.24,.72){$P''$}
  \end{picture}
\end{figure}

Consider first the case where $P$ is an acute, isosceles triangle with
apex angle of the form $\frac{\pi}{k}$, $k \in {\mathbb Z}_{\geq 3}$,
and $k$ odd.  $P$ tiles the regular $2k$--gon, $\Delta_{2k}$, with all
apex vertices at the center of $\Delta_{2k}$.  Take any copy of $P$,
call it $P'$, in this tiling.  Reflecting $P'$ across the edge opposite
the apex gives a copy, $P''$, exactly opposite $P'$ through the center of
$\Delta_{2k}$, up to translation (see Figure \ref{billiardpolyreflect}).
It follows that if we identify opposite sides of $\Delta_{2k}$ as in
Section \ref{1hangersection} we get exactly the surface and quadratic
differential (up to a complex multiple) from the billiard construction
for $P$.

Similarly, the construction from Section \ref{oddstraightchainsection}
gives the billiard surface and quadr\-atic differential for the
acute, isosceles triangle with apex angle $\frac{\pi}{k}$, $k \in
{\mathbb Z}_{\geq 3}$, and $k$ even.  When $P$ is a right triangle
with smallest angle of the form $\frac{\pi}{k}$, $k \in {\mathbb
Z}_{\geq 4}$, the construction using a regular $k$--gon in Section
\ref{1hangersection} for $k$ even, and two regular $k$--gons in Section
\ref{evenstraightchainsection} for $k$ is odd, give the billiard surface
and quadratic differential for $P$.

We have thus proved Theorem \ref{main3}, with the exception of the three
non-isosceles lattice triangles of Theorem \ref{ksptheorem}, and Theorem
\ref{main4}, except for the cases ${\cal E}_{6}$, ${\cal E}_{7}$, and
${\cal E}_{8}$.

Now, one may directly verify that $\mu(A \cup B)$ is given by
$2\cos\left(\frac{\pi}{12}\right)$, $2\cos\left(\frac{\pi}{18}\right)$,
and $2 \cos\left(\frac{\pi}{30}\right)$, for ${\cal G}(A \cup B)
= {\cal E}_{6}$, ${\cal E}_{7}$, and ${\cal E}_{8}$, respectively.
The orders of $\DAf(T_{A}T_{B})$ are thus, respectively, $6$, $9$,
and $15$.  This shows that the signatures are as required and completes
the proof of Theorem \ref{main4}.
\endproof

We now consider the billiard construction for the three
exceptional lattice triangles.  We refer to the triangles with
angles $\left(\frac{\pi}{4},\frac{\pi}{3},\frac{5 \pi}{12}\right)$,
$\left(\frac{2 \pi}{9},\frac{\pi}{3},\frac{4 \pi}{9}\right)$, and
$\left(\frac{\pi}{5},\frac{\pi}{3},\frac{7 \pi}{15}\right)$, as
$\Delta_{6}$, $\Delta_{7}$, and $\Delta_{8}$, respectively.

\begin{figure}[ht!]\small
\cl{
\psfig{file=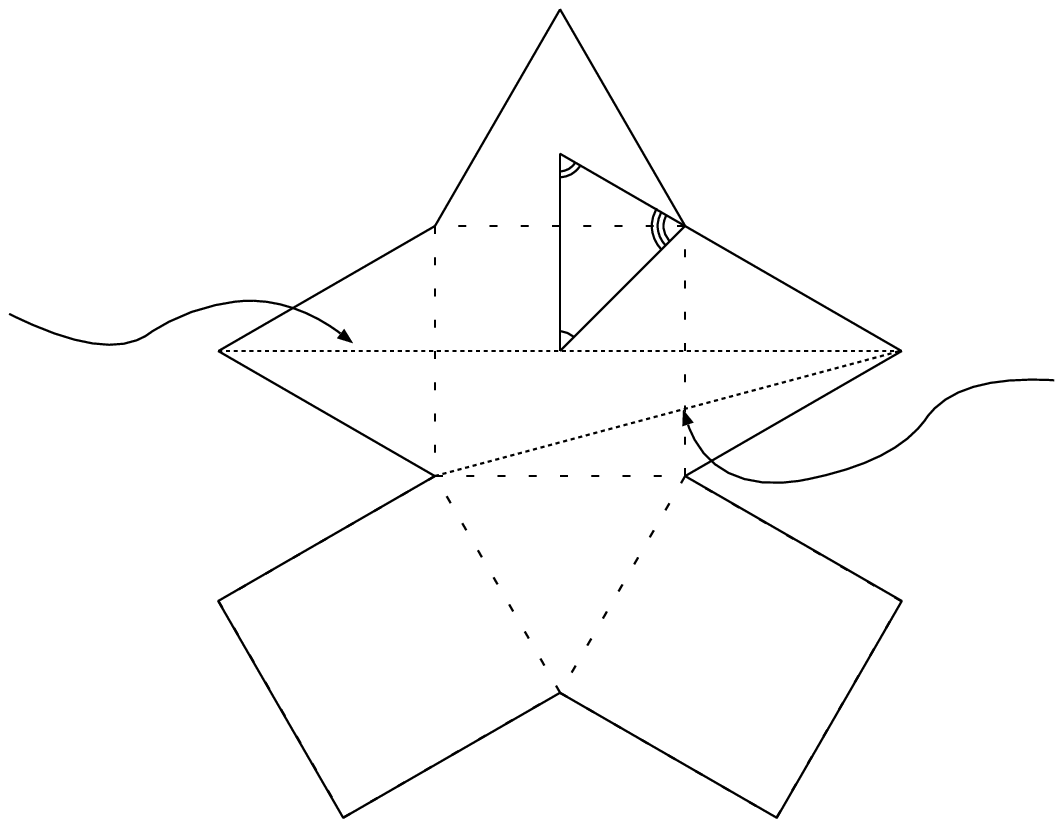,height=2.2truein}}
\caption{$\Delta_{6}$ surface}
\label{exceptional6}
\vglue 4pt
  \setlength{\unitlength}{1in}
  \begin{picture}(0,0)(0,0)
    \put(.53,1.97){$A$--direction}
    \put(4.03,1.67){$B$--direction}
    \put(1.98,2){$x$}
    \put(2.38,.63){$x$}
    \put(3.28,2.05){$y$}
    \put(2.93,.63){$y$}
    \put(2.92,2.39){$z$}
    \put(1.83,.77){$z$}
    \put(2.38,2.39){$w$}
    \put(3.46,.77){$w$}
    \put(1.98,1.34){$u$}
    \put(3.31,1.55){$u$}
    \put(3.31,1.34){$v$}
    \put(1.98,1.55){$v$}
  \end{picture}
\end{figure}

The following two surfaces and flat structures are described by Vorobets
in \cite{Vor}.  Consider first the non-convex polygon shown in Figure
\ref{exceptional6}.  This is a union of three squares and four equilateral
triangles as shown.  The surface, $S$, obtained by identifying parallel
sides as indicated has genus three, and a flat structure $q$.  $S$
is tiled by $24$ copies of $\Delta_{6}$, no two of which are parallel.
These triangles are obtained by considering the centers of the squares
and equilateral triangles and the singular point (there is just one),
and appropriately ``connecting the dots''.  One of the tiles is shown in
the Figure \ref{exceptional6}.  Since $24$ is the order of the dihedral
group generated by reflections in lines through the origin parallel
to the three sides of $\Delta_{6}$, it follows that $q$ is the flat
structure from the billiard construction for $\Delta_{6}$.

The foliations parallel to the two line segments shown in Figure
\ref{exceptional6} have all closed leaves, decomposing $S$ into annuli in
two different ways.  One direction shown is horizontal and the other makes
an angle $\frac{\pi}{12}$ with the first.  Appealing to some trigonometry
we see that these foliations do indeed define annular decompositions.
Moreover, for each of the two annular decomposition, the product of
a single Dehn twist in each annulus acts as an affine transformation.
Let $A$ denote the union of the cores of the horizontal annuli and $B$
the union of the cores of the other annuli.  One can then verify that
${\cal G}(A \cup B) = {\cal E}_{6}$.

\begin{figure}[ht!]\small
\cl{
\psfig{file=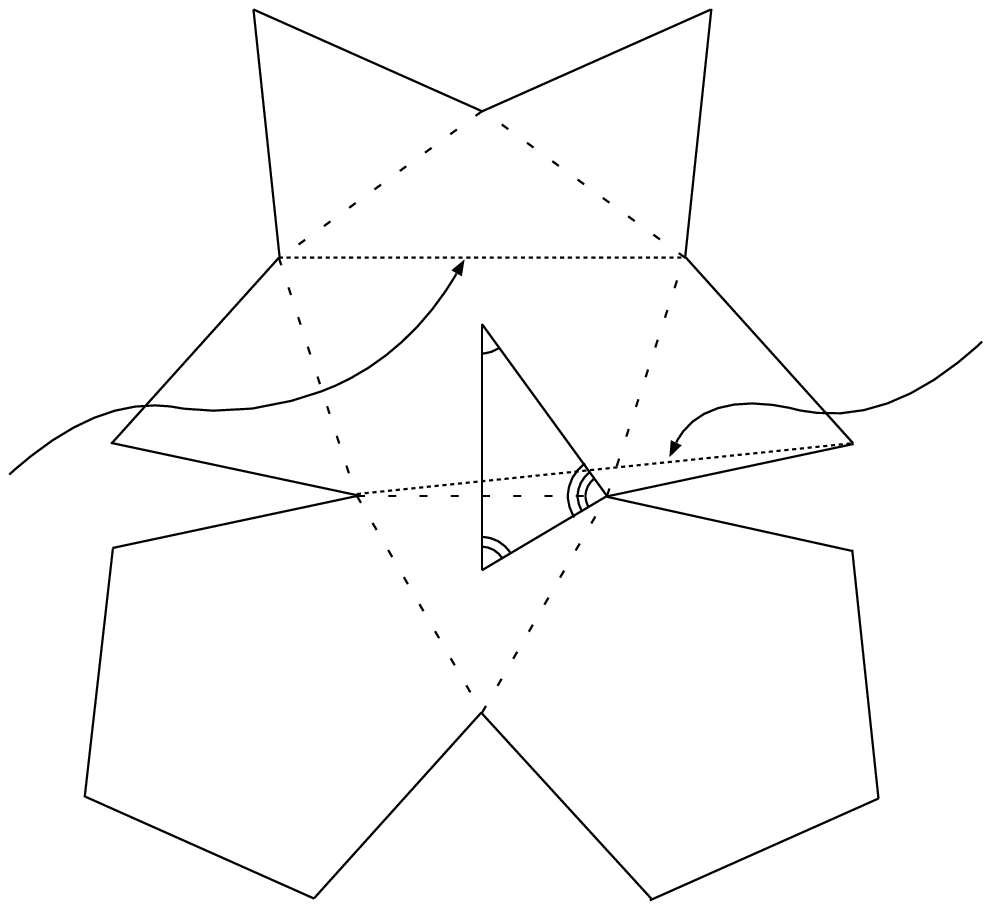,height=2.2truein}}
\caption{$\Delta_{8}$ surface}
\label{exceptional8}
\vglue 4pt
  \setlength{\unitlength}{1in}
  \begin{picture}(0,0)(0,0)
    \put(.73,1.52){$A$--direction}
    \put(3.83,1.92){$B$--direction}
    \put(2.29,2.64){$x$}
    \put(1.81,.59){$x$}
    \put(2.81,2.65){$y$}
    \put(3.28,.59){$y$}
    \put(1.78,1.89){$z$}
    \put(2.39,.71){$z$}
    \put(3.33,1.89){$w$}
    \put(2.68,.71){$w$}
    \put(1.91,1.37){$u$}
    \put(3.23,1.53){$u$}
    \put(1.87,1.53){$v$}
    \put(3.18,1.37){$v$}
    \put(1.96,2.37){$r$}
    \put(3.55,1.07){$r$}
    \put(3.13,2.37){$s$}
    \put(1.55,1.07){$s$}
  \end{picture}
\end{figure}

In a completely analogous fashion, we can consider the non-convex polygon
shown in Figure \ref{exceptional8} which is a union of three regular
pentagons and five equilateral triangles.  Identifying parallel sides as
indicated gives a surface, $S$, of genus four with a flat structure, $q$.
This is tiled by $30$, pairwise nonparallel copies of $\Delta_{8}$,
again obtained by appropriately joining centers of pentagons, regular
triangles and the singular point.  As above, we see that $q$ is the flat
structure coming from the billiard construction for $\Delta_{8}$.

One of the directions shown is horizontal and the other makes an angle
$\frac{\pi}{30}$ with the horizontal.  Again some trigonometry shows that
we obtain annular decompositions in these directions having cores $A$ and
$B$, respectively, and $T_{A}$ and $T_{B}$ act as affine transformations.
In this case ${\cal G}(A \cup B) = {\cal E}_{8}$.

\begin{figure}[ht!]\small
\cl{
\psfig{file=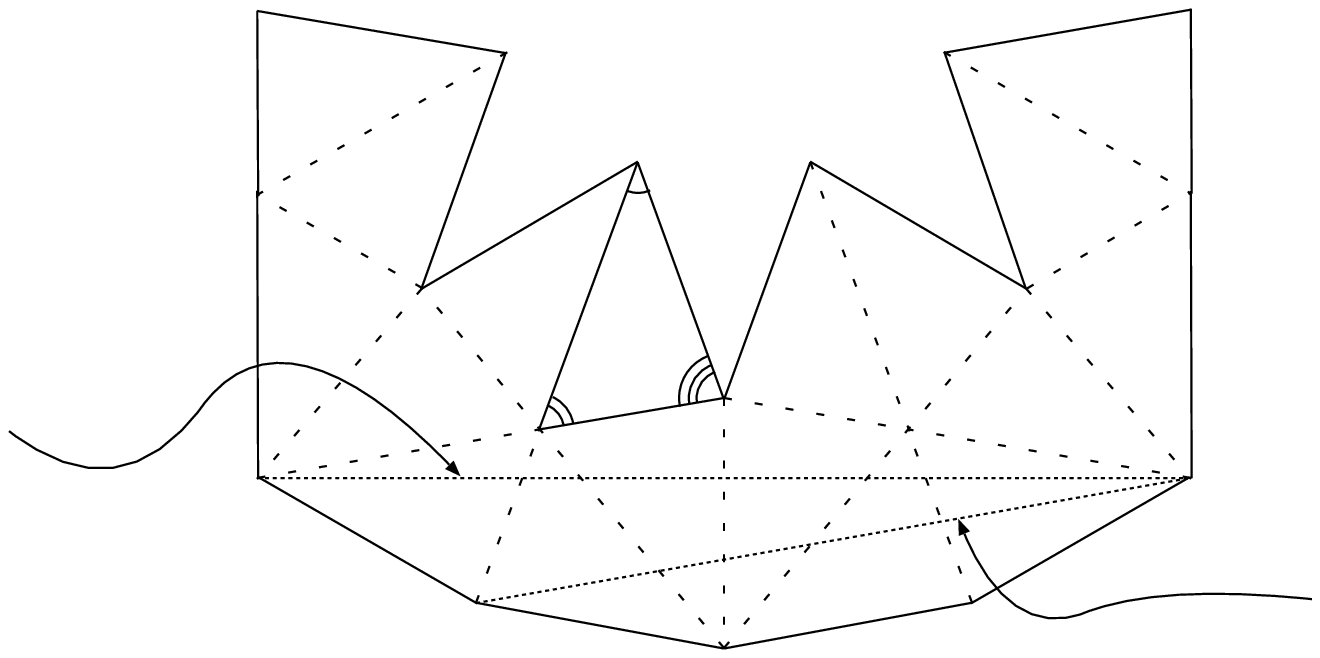,height=1.4truein}}
\caption{$\Delta_{7}$ surface}
\label{exceptional7}
\vglue 4pt
  \setlength{\unitlength}{1in}
  \begin{picture}(0,0)(0,0)
    \put(.46,1){$A$--direction}
    \put(4.06,.6){$B$--direction}
    \put(1.97,1.91){$x$}
    \put(2.4,.5){$x$}  
    \put(3.44,1.94){$y$}
    \put(3.03,.5){$y$}  
    \put(3.13,1.49){$z$}
    \put(1.89,.7){$z$}
    \put(2.23,1.49){$w$}
    \put(3.49,.71){$w$}
    \put(3.18,1.61){$u$}
    \put(2.63,1.42){$u$}
    \put(2.2,1.61){$v$} 
    \put(2.76,1.42){$v$}
    \put(1.61,1.42){$r$}
    \put(3.78,1.42){$r$}
  \end{picture}
\end{figure}

Finally, for the triangle $\Delta_{7}$, we briefly describe the
construction of Kenyon and Smillie in \cite{KS}.  We identify sides of
the polygon in Figure \ref{exceptional7} as indicated.  The result is a
surface, $S$, of genus three and a flat structure $q$.  $S$ is tiled by
$18$ copies of $\Delta_{7}$ is as indicated and again $q$ is the flat
structure from the billiard construction for $\Delta_{7}$.

The foliations in the two directions indicated (one horizontal, the
other at an angle $\frac{\pi}{18}$ from horizontal) define annular
decompositions.  For each annular decomposition, the product of a
single Dehn twist in each annulus acts as an affine transformation.
Denoting the union of the cores of the horizontal annuli by $A$ and the
union of the other cores by $B$, one can check that ${\cal G}(A \cup B)
= {\cal E}_{7}$.

Therefore, the three exceptional lattice triangles $\Delta_{6}$,
$\Delta_{7}$, and $\Delta_{8}$ define the same Teichm\"uller curves
as the configurations with graph ${\cal E}_{6}$, ${\cal E}_{7}$, and
${\cal E}_{8}$, respectively.  This completes the proof of Theorem
\ref{main3}.
\endproof

We note that for each of the realizations of recessive configurations on
surfaces we have described over the last four sections, the corresponding
quadratic differentials are squares of holomorphic $1$--forms.  Therefore,
the following is a consequence of Proposition \ref{recessiveuniqueness}
and Proposition \ref{homologyreptheorem}.

\begin{corollary}
If ${\cal G}(A \cup B)$ is recessive, then the action of $\langle
T_{A},T_{B} \rangle$ on homology is faithful.
\end{corollary}

\section{Coxeter and Artin groups} \label{coxartsection}

In this section we recall a few facts about Coxeter groups and Artin
groups which indicate a connection with groups generated by two positive
multi-twists.  We then state McMullen's Theorem \ref{mctheorem} and
verify the following.

\begin{theorem} \label{main6}
Let ${\cal G}(A \cup B)$ be non-critical dominant with small type.
Then $\sigma_{A}\sigma_{B}$ is sent by $\Psi$ to a pseudo-Anosov with
dilatation equal to the spectral radius of its image under $\Theta
\circ \pi_{ac}$.
Moreover, among all essential elements in $\langle \sigma_{A},\sigma_{B}
\rangle$, $\sigma_{A}\sigma_{B}$ minimizes both dilatation as well as
spectral radius for the respective homomorphisms.
\end{theorem}

We next examine Hironaka's Theorem \ref{hirtheorem} and use her ideas,
along with Theorem \ref{howlettthrm} of Howlett, to prove the following:

\begin{theorem} \label{main7}
Let ${\cal G}(A \cup B)$ have small type and suppose that $A$ and $B$
can be oriented so that all intersections of $A$ with $B$ are positive.
Then there exists a homomorphism
$$\eta\co {\mathbb R}^{K} \rightarrow H_{1}(S;{\mathbb R})$$
such that
$$(T_{A}T_{B})_{*} \circ \eta = - \eta \circ \Theta(\sigma_{A}\sigma_{B})$$
Moreover, $\Theta(\sigma_{A}\sigma_{B})|_{\ker(\eta)} = -I$ and $\eta$
preserves spectral radii.
\end{theorem}

This is a strengthening of the special case of Theorem \ref{main6}
in which $A$ and $B$ can be oriented as in the theorem.  For, in
this situation the flat structure defined by $A \cup B$ in Section
\ref{curves2flatsection} has no holonomy, and the dilatation of a
pseudo-Anosov is equal to the spectral radius of the action on homology
(see McMullen \cite{Mc3}).

\medskip
{\bf Remark}\qua
I would like to thank Walter Neumann for first indicating the connection
with Coxeter groups which led to Proposition \ref{genegrouptype}.

\subsection{Graphs and groups} \label{coxbasics}

Let ${\cal G}$ be any finite graph without loops (cycles of length
one), which we refer to as a {\em Coxeter graph}.  Let $\Sigma =
\{s_{1},\ldots,s_{K} \}$ be the vertices of ${\cal G}$ (throughout
Section \ref{coxartsection}, $K$ will denote the number of vertices of
${\cal G}$).  For each $1 \leq i < j \leq K$, let $m_{ij}$ be $2$ plus
the number of edges connecting $s_{i}$ to $s_{j}$, and set $m_{ii}
= 1$.  So, when ${\cal G}$ has small type, $m_{ij} \in \{ 1,2,3 \}$
for all $i,j$.  This will be the primary case of interest for us.
We will restrict ourselves to the case that ${\cal G}$ is connected
(see Humphreys \cite{Hum} for more details).

\medskip
{\bf Remark}\qua
We note that one usually allows the possibility that some vertices are
connected by infinitely many edges, but because we are only interested
in the small type case, we have not bothered to include this in the
discussion.  Also, a common convention is to consider Coxeter graphs as
graphs without multiple edges, where the edge between $s_{i}$ and $s_{j}$
is labeled with $m_{ij} \in {\mathbb Z}_{\geq 3} \cup \{ \infty \}$.
The convention we have adopted is more suitable to our situation.

Given a Coxeter graph, ${\cal G}$, there are two groups associated to it:
the {\em Coxeter group}
$${\frak C}({\cal G}) = \left\langle s_{i} \in \Sigma \, \left| \,
  (s_{i}s_{j})^{m_{ij}}\right. = 1, 1 \leq i \leq j \leq K \right\rangle$$
and the {\em Artin group}
$${\frak A}({\cal G}) = \left\langle s_{i} \in \Sigma \, \left| \,
  (s_{i}s_{j})^{\frac{m_{ij}}{2}} = (s_{j}s_{i})^{\frac{m_{ij}}{2}}\right.,
  1 \leq i < j \leq K \right\rangle$$
where for $m$ odd we define $(xy)^{\frac{m}{2}} = (xy)^{\frac{m-1}{2}}x$
(eg if $m_{ij} = 3$, the relation is the braid relation $s_{i}s_{j}s_{i}
= s_{j}s_{i}s_{j}$).

We will discuss both groups, and certain definitions are the same
in each.  However, when we wish to distinguish between an element of
${\frak A}({\cal G})$ and ${\frak C}({\cal G})$, we will denote the
former with a ``prime''.  Thus, $s_{i} \in {\frak C}({\cal G})$ and
$s_{i}' \in {\frak A}({\cal G})$.

${\frak C}({\cal G})$ is obtained from ${\frak A}({\cal G})$ by adding
the relation $s_{i}^{2} = 1$.  This defines an epimorphism
$$\pi_{ac}\co {\frak A}({\cal G}) \rightarrow {\frak C}({\cal G})$$
obtained by sending $s_{i}'$ to $s_{i}$.

Given a Coxeter or Artin group a {\em special subgroup} is any subgroup
generated by a subset $\Sigma_{0} \subset \Sigma$.  These special
subgroups are precisely the Coxeter and Artin groups, respectively,
associated to the largest subgraph of ${\cal G}$ having $\Sigma_{0}$
as its vertex set.  An element of ${\frak C}({\cal G})$ is said to be
{\em essential} if it is not conjugate into any proper special subgroup,
and we call an element ${\frak A}({\cal G})$ essential if its image by
$\pi_{ac}$ is essential.

The product of all the generators (in any order) is called a {\em Coxeter
element}.  For any bipartite graph ${\cal G}$, there is a special Coxeter
element defined as follows.  Since the graph is bipartite there exists a
partition $\Sigma = A \cup B$ so that no two $A$--vertices (respectively,
$B$--vertices) are adjacent.  The product of the elements of $A$
(respectively, $B$) defines
$$\sigma_{A} = \prod_{s_{i} \in A} s_{i} \quad \mbox{ and } \quad
  \sigma_{B} = \prod_{s_{j} \in B} s_{j}.$$
The product $\sigma_{A}\sigma_{B}$ is called the {\em bi-colored Coxeter
element}.

\subsection{Artin groups and mapping class groups} \label{artmapsection}

Let $A,B \in {\cal S}'(S)$ with ${\cal G} = {\cal G}(A \cup B)$ of
small type.  There is a nice relationship between ${\frak A}({\cal G})$
and $\Mod(S)$.  The vertices $s_{1},\ldots,s_{K}$ of ${\cal G}$ can be
identified with the components of $A$ and $B$ as well as generators
of ${\frak A}({\cal G})$.  By (\ref{braidrelationequation}) from Section
\ref{multitwistsection} and the definition of ${\frak A}({\cal G})$,
we can define a homomorphism
$$\Psi\co {\frak A}({\cal G}) \rightarrow \Mod(S)$$
by sending the generator $s_{i}$ to the Dehn twist about the curve
corresponding to the vertex $s_{i}$.  Note that, after relabeling if
necessary, we have $\Psi(\sigma_{A}) = T_{A}$ and $\Psi(\sigma_{B})
= T_{B}$.

\medskip
{\bf Remark}\qua This construction can be carried out for any graph of
small type (not necessarily bipartite).
Indeed, to such a graph one can (nonuniquely) associate a surface
and a set of curves, pairwise intersecting at most once, and define
a homomorphisms $\Psi$ as above.  For more on this see eg Crisp and
Paris \cite{CP}, Wajnryb \cite{Waj1}, \cite{Waj2}, Perron and Vannier
\cite{PV}, and A'Campo \cite{AC4}.

\subsection{Geometric representations of Coxeter groups}
\label{geomrepcoxsection}

Suppose that ${\cal G}$ is a connected Coxeter graph.  There is an
associated quadratic form $\Pi_{{\cal G}}$ on ${\mathbb R}^{K}$ and a
faithful representation
$$\Theta\co{\frak C}({\cal G}) \rightarrow \On(\Pi_{{\cal G}})$$
where $\On(\Pi_{{\cal G}})$ is the orthogonal group of the quadratic form
$\Pi_{{\cal G}}$, and each generator $s_{i} \in \Sigma$ is represented
by a reflection.

Up to equivalence over ${\mathbb R}$, there are precisely four possibilities for the form $\Pi_{{\cal G}}$ (see \cite{Hum}).
These are characterized by the signature, $sgn(\Pi_{{\cal G}})$, and $K$.
Accordingly, the group ${\frak C}({\cal G})$ is said to be:
\begin{itemize}
\item {\em spherical} if $sgn(\Pi_{{\cal G}}) = (K,0)$,
\item {\em affine} if $sgn(\Pi_{{\cal G}}) = (K-1,0)$,
\item {\em hyperbolic} if $sgn(\Pi_{{\cal G}}) = (p,1)$ and $p + 1 \leq K$, and
\item {\em higher-rank} if $sgn(\Pi_{{\cal G}}) = (p,q)$ and $p + q \leq K$, $q \geq 2$.
\end{itemize}

When ${\cal G}$ has small type (our only case of interest), this quadratic
form is easily described in terms of ${\cal A}d({\cal G})$, the associated
adjacency matrix (see Section \ref{graphsubsection}).
The form $\Pi_{\cal G}$ is then defined by the matrix
$$2I - {\cal A}d({\cal G}).$$
Because ${\cal G}$ is connected, ${\cal A}d({\cal G})$ is an irreducible
matrix.  Moreover, one can easily see that the group ${\frak C}({\cal G})$
is spherical or affine if and only if $\mu({\cal G}) < 2$ or $\mu({\cal
G}) = 2$, respectively.  The following is therefore a consequence of
Theorem \ref{smithlist}.

\begin{proposition} \label{genegrouptype}
For ${\cal G} = {\cal G}(A \cup B)$ of small type we have
\begin{itemize}
\item ${\cal G}$ is recessive if and only if ${\frak C}({\cal G})$ is spherical,
\item ${\cal G}$ is critical if and only if ${\frak C}({\cal G})$ is affine, and
\item ${\cal G}$ is non-critical dominant if and only if ${\frak C}({\cal G})$ is hyperbolic or higher-rank.
\end{itemize}
\end{proposition}

This proposition and the construction mentioned in the previous section
begin to shed light on an interesting connection between Coxeter and
Artin groups and the work presented so far.  The following result of
McMullen \cite{Mc} indicates that the connection is stronger still
(see also \cite{BLM} and \cite{AC3}).

\begin{theorem}[McMullen]
\label{mctheorem}
Suppose ${\cal G}$ is bipartite and has small type.
Then over all essential $\phi \in {\frak C}({\cal G})$, the spectral
radius of $\Theta(\phi)$ is minimized by the bi-colored Coxeter element
$\phi = \sigma_{A}\sigma_{B}$, and is given as the larger absolute value
of a root of the polynomial
$$x^{2} + x(2 - (\mu({\cal G}))^{2}) + 1$$
For spherical or affine Coxeter groups, this minimum is $1$, and among
all hyperbolic or higher-rank Coxeter groups, the minimal spectral radius
is uniquely minimized for the Coxeter group ${\frak C}({\cal E}h_{10})$,
and the minimal spectral radius is precisely $\lambda_{L}$.
\end{theorem}

Theorem \ref{main6} is an immediate consequence of Theorem \ref{main5},
Proposition \ref{genegrouptype}, and Theorem \ref{mctheorem}.
\endproof

{\bf Remarks}

(1)\qua McMullen's theorem does not require ${\cal G}$ to have small type,
although in that case the appearance of $\mu({\cal G})$ in the theorem
is replaced by the spectral radius of the {\em Coxeter adjacency matrix}.
The bipartite assumption is also unnecessary, however in this case, the
number given by the theorem is only a lower bound, not the minimum. 

(2)\qua We note that for a random element of $\langle \sigma_{A},\sigma_{B}
\rangle < {\frak A}({\cal G})$, there is no connection between is spectral
radius under $\Theta \circ \pi_{ac}$ and its dilatation under $\Psi$.
In particular, solving the minimization problem for one does not solve
it for the other.

\subsection{Coxeter links and Coxeter actions}  \label{coxlinksection}

In \cite{Hir1}, Hironaka provides a construction of a fibered link in
$S^{3}$ which depends on a Coxeter graph of small type (as well as some
additional data).
This expresses the link complement as the mapping torus of an automorphism
of the fiber called the {\em monodromy}.
The main theorem of \cite{Hir1} states that, up to sign, the action
on homology of the monodromy is conjugate to the geometric action of a
certain Coxeter element (see below for the precise statement).
As we shall see, under certain additional hypotheses, the monodromy is
of the form $T_{A}T_{B}$ for appropriate $A \cup B$ filling the fiber.

\begin{figure}[ht!]\small\vglue 6pt
\cl{
\psfig{file=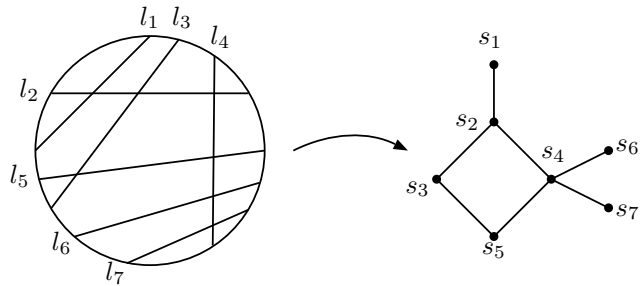,height=1.25truein}}
\caption{A Coxeter graph from a chord diagram}
\label{chorddi1}
\vglue 4pt
  \setlength{\unitlength}{1in}
  \begin{picture}(0,0)(0,0)
    \put(1.03,1.44){$l_{2}$}
    \put(1.66,1.82){$l_{1}$}
    \put(1.83,1.82){$l_{3}$}
    \put(2.03,1.72){$l_{4}$}
    \put(.98,.99){$l_{5}$}
    \put(1.2,.64){$l_{6}$}
    \put(1.48,.48){$l_{7}$}
    \put(3.43,1.72){$s_{1}$}
    \put(3.31,1.29){$s_{2}$}
    \put(3.05,.94){$s_{3}$}
    \put(3.76,1.12){$s_{4}$}
    \put(3.45,.62){$s_{5}$}
    \put(4.15,1.16){$s_{6}$}
    \put(4.15,.82){$s_{7}$}
  \end{picture}
\end{figure}

We now describe Hironaka's construction (for more details, see
\cite{Hir1}).
A {\em chord diagram} is a collection of straight arcs ${\cal L} = \{
l_{1},\ldots,l_{K} \}$, called {\em chords}, in the unit disk ${\mathbb D}
\subset {\mathbb R}^{2}$ connecting mutually disjoint pairs of points
on the boundary of ${\mathbb D}$.
The chord diagram defines a Coxeter graph ${\cal G}$ of small type
as follows.
The vertices $\Sigma$ are identified with the chords of ${\cal L}$, and
two vertices $s_{i}$ and $s_{j}$ are joined by an edge if and only if the
corresponding chords intersect non-trivially (see Figure \ref{chorddi1}).
For a chord diagram ${\cal L}$ defining a Coxeter graph ${\cal G}$,
we say that ${\cal L}$ is a {\em chord realization} of ${\cal G}$,
and that the graph is {\em chord-realizable}.

\begin{figure}[ht!]\small
\cl{
\psfig{file=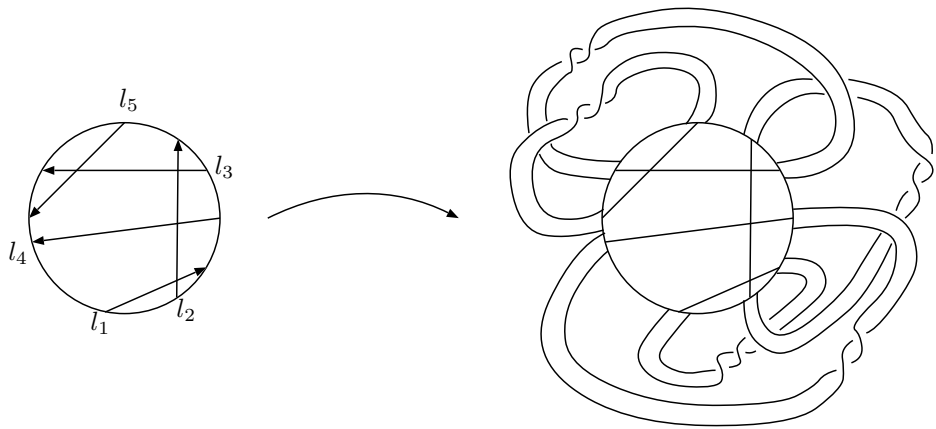,height=2.25truein}}
\caption{Plumbing Hopf bands onto a chord diagram}
\label{chorddi2}
\vglue4pt
  \setlength{\unitlength}{1in}
  \begin{picture}(0,0)(0,0)
    \put(1.03,1.12){$l_{2}$}
    \put(1.23,1.87){$l_{3}$}
    \put(.73,2.22){$l_{5}$}
    \put(.15,1.42){$l_{4}$}
    \put(.58,1.05){$l_{1}$}
  \end{picture}
\end{figure}

Suppose that ${\cal G}$ is a chord-realizable Coxeter graph.  An ordering
on the vertices $\Sigma = \{ s_{1},\ldots,s_{K} \}$ (equivalently, an
ordering on the chords ${\cal L} = \{ l_{1},\ldots, l_{K} \}$) gives
rise to a fibered link as follows.  Recall that a {\em Hopf band} $H$
is an annulus spanning a Hopf link $L$ (see Figure \ref{chorddi2}).
For each chord we plumb a right-handed Hopf band onto the disk in
$S^{3}$ so that the core of the band agrees with the chord in the disk.
We do this in the order specified by the ordering of the chords (see
Figure \ref{chorddi2}).  We denote the resulting surface by $S$, and its
boundary by $L = \partial S$.  It is well known that $L$ is a fibered
link (see \cite{Gab}).  We also note that the ordering of the vertices
also specifies a Coxeter element $c = s_{1}s_{2} \cdots s_{K}$.

Finally, if we orient the chords in a chord diagram so that the ordering
is {\em compatible} with the orientation, then the resulting link is
said to be a {\em Coxeter link}.  The compatibility here simply means
that for $i < j$, the chord $s_{i}$ must intersect the chord $s_{j}$
positively (if at all).  For example, the ordering of the chords in
Figure \ref{chorddi2} is compatible with the orientations.

The following is proved in \cite{Hir1}.

\begin{theorem}[Hironaka] \label{hirtheorem}
Given an oriented, ordered chord diagram with associated Coxeter graph
${\cal G}$, Coxeter link $L = \partial S$, fiber $S$, and monodromy
$\phi$, there exists an isomorphism
$$\nu\co{\mathbb R}^{K} \rightarrow H_{1}(S;{\mathbb R})$$
such that $\phi_{*} \circ \nu = - \nu \circ \Theta(c)$
where $c$ is the Coxeter element determined by the ordering.
If the spectral radius of $\phi_{*}$ is greater than $1$, then it is
bounded below by $\lambda_{L}$.
\end{theorem}

Hironaka's proof uses the following interpretation of a theorem of Howlett \cite{How} in the case that ${\cal G}$ has small type.

\begin{theorem}[Howlett]
\label{howlettthrm}
If ${\cal G}$ is of small type, then the Coxeter element $c$ is given by
$$c = -(I - {\cal A}d({\cal G})^{+})^{-1}(I - {\cal A}d({\cal G})^{+})^{t}$$
where ${\cal A}d({\cal G})^{+}$ is the upper triangular part of the
adjacency matrix ${\cal A}d({\cal G})$.
\end{theorem}

To prove her theorem, Hironaka shows that the Seifert matrix for the
link is given by $I - {\cal A}d({\cal G})^{+}$.  It then follows from
classical knot theory (see eg \cite{R}) that the action of the monodromy
on homology is given by $c$, as required.

\medskip
{\bf Remark}\qua There is another construction of fibered links for
which Dehn twists and Coxeter diagrams appear very naturally.
This is described by A'Campo in \cite{AC1}, \cite{AC2}, and the references
contained therein.
Although we have not fully investigated this, it seems likely that this
construction is closely related to the one described above.

To relate Hironaka's Theorem to our work, we recall that according
to Gabai \cite{Gab2} the monodromy of a fibered link obtained by
(generalized) plumbing of two fibers is the composition of the two
monodromies (see \cite{Gab2} for a more precise statement).  The monodromy
for a Hopf link (with fiber a right-handed Hopf band) is a positive Dehn
twist about the core of the band.  It follows that the monodromy of the
Coxeter link constructed above is the product of Dehn twists about the
cores of the plumbed on Hopf bands (the product is taken, from left to
right, in the order given by the ordering of the chords).

Suppose now that a chord diagram has bipartite Coxeter graph ${\cal G}$
with vertices $\Sigma$ colored by $A$ and $B$, and there is an ordering
of the vertices so that for all $s_{i} \in A$ and $s_{j} \in B$, we have
$i < j$.  We call this a {\em bi-colored ordering} with respect to $A$
and $B$.  The cores of the Hopf bands associated to $A$ (respectively,
$B$) give an essential $1$--manifold we also denote by $A$ (respectively,
$B$) in the surface $S$.  It is easy to see that ${\cal G} = {\cal G}(A
\cup B)$.

The previous two paragraphs imply the following theorem.

\begin{theorem} \label{bicoloredplumber}
In the setting of Theorem \ref{hirtheorem}, if we further assume that
the ordering is a bi-colored ordering with respect to $A$ and $B$,
then $\phi = T_{A}T_{B}$.
In particular, the action of $(T_{A}T_{B})_{*}$ on $H_{1}(S;{\mathbb
R})$ is conjugate to the action of $-\Theta(\sigma_{A}\sigma_{B})$
on ${\mathbb R}^{K}$.
\end{theorem}

It is not hard to see that this implies Theorem \ref{main7} in the
special case that $T_{A}T_{B}$ is the monodromy for a Coxeter link.
We now give the proof in the general case.

\medskip
{\bf Proof of Theorem \ref{main7}}\qua
Suppose ${\cal G}(A \cup B)$ has small type and we have oriented $A$
and $B$ so that all intersections of $A$ with $B$ are positive.
So, for any pair of components $a_{i} \subset A$ and $b_{j} \subset B$, we have
\begin{equation} \label{alggeominteqn}
a_{i} \cdot b_{j} = i(a_{i},b_{j}) = -b_{j} \cdot a_{i}.
\end{equation}
Let ${\cal N}(A \cup B)$ be a regular neighborhood of $A \cup B$ in
${\cal S}$, and let us denote the inclusion into $S$ by $\iota\co{\cal
N}(A \cup B) \rightarrow S$.

We may define a monomorphism
$$\eta_{0}\co {\mathbb R}^{K} \rightarrow H_{1}({\cal N}(A
  \cup B);{\mathbb R})$$
by sending each basis element to the homology class of the oriented
curve it determines.
We view $a_{1},\ldots,a_{n},b_{1},\ldots,b_{m}$ as an ordered basis for both
${\mathbb R}^{K}$ as well as the subspace $V \subset H_{1}({\cal N}(A
\cup B);{\mathbb R})$ which they span.
To see that these are indeed linearly independent in $H_{1}({\cal N}(A
\cup B);{\mathbb R})$, we note that for any one of these, say $a_{1}$,
we can easily find an arc $\alpha$ which intersects $a_{1}$ once, but
misses all the others.
This arc determines an element of $H_{1}({\cal N}(A \cup B),\partial
{\cal N}(A \cup B); {\mathbb R})$, which by Poincar\'e duality, is
identified with the dual space of $H_{1}({\cal N}(A \cup B);{\mathbb R})$
via intersection numbers.
Thus there is a functional vanishing on all the vectors except $a_{1}$.
Since $a_{1}$ was arbitrary, the vectors are linearly independent.

The action on homology of a Dehn twist $T_{a}$ is given by
$$(T_{a})_{*}(x) = x + (a \cdot x) a.$$
By (\ref{alggeominteqn}), the matrix for the actions of $T_{A}$ and
$T_{B}$ on $V$ with respect to the basis $a_{1},\ldots,a_{n},b_{1},\ldots,b_{m}$
is thus given by
$$(T_{A})_{*} = \left( \begin{array} {cc}
  I & N \\
  0 & I \\ \end{array} \right) \quad \mbox{ and } \quad
  (T_{B})_{*} = \left( \begin{array} {cc}
  I & 0 \\
  -N^{t} & I \\ \end{array} \right)$$
where $N_{ij} = i(a_{i},b_{j})$ as in Section \ref{curves2flatsection}.

Now, by Theorem \ref{howlettthrm} and Proposition \ref{whatismu}, we have:
$$\begin{array}{ccl}
  \Theta(\sigma_{A}\sigma_{B}) &  =  & -(I-{\cal A}d({\cal G})^{+})^{-1}(I
  - {\cal A}d({\cal G})^{+})^{t} \\
   & & \\
   & = & -\left( \begin{array}{cc}
  I & -N \\
  0 & I \\ \end{array} \right)^{-1}
  \left( \begin{array}{cc}
  I & -N \\
  0 & I \\ \end{array} \right)^{t}
  = - \left( \begin{array}{cc}
  I & N \\
  0 & I \\ \end{array} \right)
  \left( \begin{array}{cc}
  I & 0 \\
  -N^{t} & I \\ \end{array} \right)\\
   & & \\
   & = & -(T_{A})_{*} (T_{B})_{*} = -(T_{A}T_{B})_{*}\\ \end{array}$$
The matrices are given with respect to $a_{1},\ldots,a_{n},b_{1},\ldots,b_{m}$
in both ${\mathbb R}^{K}$ and $H_{1}({\cal N}(A \cup B);{\mathbb R})$.
Because $\eta_{0}$ is ``the identity'' with respect to this basis,
we see that
\begin{equation} \label{eta0eqn}
(T_{A}T_{B})_{*} \circ \eta_{0} = -\eta_{0} \circ \Theta(\sigma_{A}\sigma_{B}).
\end{equation}
We now obtain the required homomorphism
$$\eta = \iota_{*} \circ \eta_{0}\co{\mathbb R}^{K} \rightarrow
  H_{1}(S;{\mathbb R}).$$
Because $T_{A}T_{B}$ is supported on ${\cal N}(A \cup B)$, we have $\iota
\circ T_{A}T_{B} = T_{A}T_{B} \circ \iota$, and hence
$$(T_{A}T_{B})_{*} \circ \eta = -\eta \circ \Theta(\sigma_{A}\sigma_{B}).$$
This proves Theorem \ref{main7}, except for the last sentence.

To see this, note that the kernel of $\eta$ is the image of the kernel
of $\iota_{*}$ by $\eta_{0}^{-1}$.  Since $S$ is obtained from ${\cal
N}(A \cup B)$ by gluing disks to the boundary, we see that the kernel of
$\iota_{*}$ consists of the span of the homology classes of the boundary.
However, $T_{A}T_{B}$ fixes the boundary pointwise, and so acts as
the identity on this span in $H_{1}$.  Therefore, by (\ref{eta0eqn}),
$\Theta(\sigma_{A}\sigma_{B})$ acts as $-I$ on the kernel of $\eta$.

Finally, we see that spectral radii are preserved in the pseudo-Anosov
case by Theorem \ref{main6} and the fact that the dilatation is equal
to the spectral radius (see \cite{Mc3}).  The only other case is when
${\cal G}(A \cup B)$ is recessive or critical.  In these cases, the
spectral radius (which is $1$) is necessarily preserved.
\endproof

\section{Applications and questions} \label{applicationsection}

Here we provide a few applications of our work and state a few interesting
questions.

\subsection{Lehmer's question, Salem numbers, and Teichm\"uller curves}
\label{lehmersapplication}

The interest in Lehmer's number stems from a problem in number theory
known as Lehmer's question (see \cite{Leh}).  To state it, we recall
that given a monic integral polynomial $p(x) \in {\mathbb Z}[x]$, the
{\em Mahler measure} of $p$ is defined by
$$\Omega(p) = \prod_{p(\theta) = 0} \max\{1,|\theta|\}.$$

\begin{question}[Lehmer]
\label{lehmersconjecture}
Is there an $\epsilon > 1$ such that if $\Omega(p) > 1$, then $\Omega(p)$\break
$\geq \epsilon$?
\end{question}

At present, the smallest known Mahler measure greater than
$1$ occurs for Lehmer's polynomial (\ref{lehmerseqn}), Section
\ref{introlehmercoxsection}, and is equal to $\lambda_{L}$.  One may view
Theorem \ref{main5} as a resolution of Lehmer's question in a particular
situation.  More precisely, if we let ${\frak D}_2$ denote the set of
all dilatations of pseudo-Anosov elements in groups generated by two
positive multi-twists, then Theorem \ref{main5} implies the following.

\begin{corollary}
The Mahler measure of the minimal polynomial of any element of ${\frak
D}_{2}$ is bounded below by $\lambda_{L}$.
\end{corollary}

A {\em Salem number} is an algebraic integer $\lambda > 1$, such that
the Galois conjugates include $\lambda^{-1}$ and all (except $\lambda$)
lie in the unit disk.  Note that a Salem number is equal to the Mahler
measure of its minimal polynomial.  In particular, an affirmative answer
to the following (see \cite{Boy1} and \cite{GH}) would be a consequence
of such an answer to Lehmer's question.

\begin{question}
\label{salemsconjecture}
Is there an $\epsilon > 1$ such every Salem number $\lambda$ satisfies
$\lambda \geq \epsilon$?
\end{question}

Lehmer's number is a Salem number, so of course the best guess
for $\epsilon$ is $\lambda_{L}$.  Because of this question, one is
generally interested in ``small'' Salem numbers.  There are currently
$47$ known Salem numbers less than $1.3$, including $\lambda_{L}$ (see
\cite{Boy1}, \cite{Boy2}, \cite{Mos}, and also \cite{FGR}).  However,
we only obtain $5$ small Salem numbers as elements of ${\frak D}_{2}$.
This set consists of all but $1$ of the Salem numbers obtained
by McMullen in \cite{Mc} as spectral radii of certain elements of
Coxeter groups.  This is not surprising, given Theorem \ref{main5}
and the fact that $5$ of the $6$ small Salem numbers obtained by
McMullen come from bicolored Coxeter elements.  On the other hand,
the only dilatations in ${\frak D}_{2}$ which can occur in the interval
$(1,1.3)$ are of the form $\lambda(T_{A}T_{B})$, for ${\cal G}(A \cup
B)$ of small type (see the proof of Theorem \ref{cdgbnmin} and use
Proposition \ref{shortesttrianglegeodesics} along with the fact that
there is exactly one simple closed geodesic when ${\cal G}(A \cup B)$
is non-critical dominant).

The elements of ${\frak D}_{2}$ are not at all representative of the
general case of dilatations of pseudo-Anosov automorphisms which are
{\em not} bounded away from $1$ (see Penner \cite{P}, Bauer \cite{Ba},
and McMullen \cite{Mc2}).  In particular, we ask the following:

\begin{question}
Which Salem numbers occur as dilatations of pseudo-Anosov automorphisms?
\end{question}

\begin{question}
Is there some topological condition on a pseudo-Anosov which guarantees
that its dilatation is a Salem number?
\end{question}

In the same vein as Questions \ref{lehmersconjecture} and
\ref{salemsconjecture}, we ask the following:

\begin{question}
Is there an $\epsilon > 1$, such that if $\phi$ is a pseudo-Anosov
automorphism in a finite co-area Teichm\"uller disk stabilizer, then
$\lambda(\phi) \geq \epsilon$?
\end{question}

Given that the dilatations we are obtaining are naturally occurring as
spectral radii of hyperbolic elements in certain non-elementary Fuchsian
groups, we would be remiss not to mention the following (see \cite{NR},
\cite{MR}, and also \cite{GH}).

\begin{theorem}[Neumann--Reid]
\label{neumannreidtheorem}
The Salem numbers are precisely the spectral radii of hyperbolic elements
of arithmetic Fuchsian groups derived from quaternion algebras.
\end{theorem}

However, because the non-cocompact arithmetic Fuchsian groups are
necessarily commensurable with $PSL_{2}{\mathbb Z}$, relatively few
of the groups generated by two positive multi-twists even inject into
arithmetic groups.

\subsection{Unexpected multi-twists and the $3$--chain relation}
\label{3chainapplication}

The work in this paper has a connection to a problem posed by McCarthy at
the 2002 AMS meeting in Ann Arbor, MI.  This was to determine the extent
to which the lantern relation in the mapping class group is characterized
by its algebraic properties (in particular the intersection patterns of
the defining curves).  Two different solutions to this were obtained,
independently by Hamidi-Tehrani in \cite{H}, and by Margalit in \cite{M}.

This question asks us to decide when an element in a group generated by
two Dehn twists can be a multi-twist.  One could ask the same question
more generally, ie for positive multi-twists.  The answer is given, to a
certain extent, by Proposition \ref{elementtypes}, Theorem \ref{main4},
and Theorem \ref{reptheorem}.  We do not spell this out here, but
instead provide a partial answer to a related question posed by Margalit
in \cite{M}.  I am grateful to Joan Birman for pointing this out to me.

Margalit asks to what extent the {\em $n$--chain relation} can be
characterized.  This is the relation
$$(T_{a_{1}}T_{a_{2}} \cdots T_{a_{n}})^{k} = M$$
where:
\begin{itemize}
\item $a_{1},\ldots,a_{n}$ are essential simple closed curves on a surface
  with $i(a_{i},a_{i+1}) = 1$, $i=1,\ldots,n-1$, and all other intersection
  numbers $0$,
\item $M$ is either $T_{d}$ or $T_{d_{1}}T_{d_{2}}$, where $d$ or
  $d_{1} \cup d_{2}$ is the boundary of a regular neighborhood of $a_{1}
  \cup \cdots \cup a_{n}$ (depending on whether $n$ is even or odd,
respectively), and
\item $k=2n+2$ for $n$ even, and $k=n+1$ for $n$ odd.
\end{itemize}

Margalit gives a characterization for $n = 2$, which we state here.

\begin{theorem}[Margalit]
Suppose $M = (T_{x}T_{y})^{k}$, where $M$ is a multi-twist and $k \in
{\mathbb Z}$, is a non-trivial relation between powers of Dehn twists
in $\Mod(S)$, and $[M,T_{x}] = 1$.
Then the given relation is the $2$--chain relation, ie $M = T_{c}^{j}$,
where $c$ is the boundary of a neighborhood of $x \cup y$, $i(x,y) =
1$, and $k = 6j$.
\end{theorem}

We note that although our work has been primarily concerned with
groups generated by two multi-twists, we can in fact obtain a similar
characterization of the $3$--chain relation.

\begin{theorem}
Suppose $M = (T_{x}T_{y}T_{z})^{k}$, where $M$ is a multi-twist and
$k \in {\mathbb Z}$, is a non-trivial relation between powers of Dehn
twists in $\Mod(S)$, and $[M,T_{x}] = [T_{x},T_{z}] = 1$.
Then the given relation is the $3$--chain relation, ie $M =
(T_{c}T_{d})^{j}$, where $c \cup d$ is the boundary of a neighborhood
of $x \cup y \cup z$, $i(x,y) = i(y,z) = 1$, and $k = 4j$.
\end{theorem}

The non-triviality here means that $i(x,y) \neq 0$, $i(y,z) \neq 0$,
and $k \neq 0$.

\begin{proof}
Since $[T_{x},M] = 1$, conjugating by $T_{x}^{-1}$, we obtain
$$M = T_{x}^{-1}(T_{x}T_{y}T_{z})^{k}T_{x} = (T_{y}T_{z}T_{x})^{k} =
  (T_{y} (T_{z}T_{x}))^{k}.$$
Also, $T_{z}T_{x}$ is a positive multi-twist since $[T_{x},T_{z}] = 1$
implies $i(x,z) = 0$.

Proposition \ref{elementtypes} tells us that ${\cal G}(y \cup (z \cup
x))$ is recessive or critical (otherwise $T_{y}T_{z}T_{x}$, and all of
its powers, would be pseudo-Anosov on the subsurface filled by $y \cup
(z \cup x)$, contradicting the fact that some power is a multi-twist).
The only such graph with $3$ vertices is ${\cal A}_{3}$.
\end{proof}

\section{Appendix A: The nonfilling case} \label{appendixreduction}

Here we provide a proof of the following.

\begin{proposition}
\label{reductionproposition}
It suffices to prove Theorem \ref{main1} for $A \cup B$ filling $S$.
\end{proposition}

As the proof will require us to deal with surfaces having nonempty
boundary, we can also allow $S$ to have boundary with no added
complications.  In particular, Theorem \ref{main1} remains true in
this setting.  Now, an allowable homeomorphism is one which leaves the
marked points invariant and fixes the boundary components pointwise.
The definition of $\Mod(S)$ is as in Section \ref{autosection}.

On $S$, consider two elements $A,B \in {\cal S}'(S)$.  Let ${\cal
N}(A \cup B)$ denote the regular neighborhood of $A \cup B$ in $S$.
Write $\overline{S} = \overline{S}_{A \cup B}$ for the subsurface of $S$
obtained by taking the union of ${\cal N}(A \cup B)$ with any open disks,
once-marked open disks, and half open annuli in the complement of ${\cal
N}(A \cup B)$.

Next, let $\widehat{S}$ be the surface obtained from $\overline{S}$ by
gluing a disk with one marked point to each boundary component, and write
$$\varepsilon\co\overline{S} \rightarrow \widehat{S}$$
for the inclusion.
We let $A$ and $B$ denote the images under $\epsilon$ of the
$1$--manifolds of the same name.
The components of $\widehat{S}$ bijectively correspond to the components,
$A_{1} \cup B_{1},\ldots,A_{k} \cup B_{k}$, of $A \cup B$, and we write
these as $\widehat{S}_{1},\ldots,\widehat{S}_{k}$.

Note that $A_{r} \cup B_{r}$ fills each component $\widehat{S}_{r}$,
except when $\widehat{S}_{r} \cong S_{0,2}$.  In this situation
$A_{r} \cup B_{r}$ is a single closed curve which is not essential in
$\widehat{S}_{r}$.  However, it should be clear from what follows that
this technicality may be ignored.

The groups we need to consider are
$$ \begin{array} {lll}
  G = \langle T_{A},T_{B}  \rangle < \Mod(S) & \hspace{2cm} & \overline{G}
  = \langle T_{A} , T_{B}  \rangle < \Mod(\overline{S}) \\
  \widehat{G} = \langle T_{A} , T_{B}  \rangle < \Mod(\widehat{S})
  & & \widehat{G}_{r} = \langle T_{A_{r}} , T_{B_{r}}  \rangle <
  \Mod(\widehat{S}_{r})\\
  \end{array}$$
for each $r = 1,\ldots,k$.  $G$ is the group from Theorem \ref{main1}.

Proposition \ref{reductionproposition} follows easily from the next
proposition since $A_{r} \cup B_{r}$ fills each of $\widehat{S}_{r}$.

\begin{proposition}
\label{reduction2filling}
For $G$ and $\widehat{G}_{1},\ldots,\widehat{G}_{k}$ as above
$$G \cong {\mathbb F}_{2} \Leftrightarrow \widehat{G}_{r}
  \cong {\mathbb F}_{2}$$
for some $r = 1,\ldots,k$.
\end{proposition}

\begin{proof}
The map $\varepsilon$ induces an epimorphism
$$\varepsilon_{*}\co\Mod(\overline{S}) \rightarrow \Mod(\widehat{S}).$$
Moreover, the kernel of $\varepsilon_{*}$ is generated by Dehn twists
about curves parallel to the boundary components of $\overline{S}$,
which defines the following central extension \cite{Bic}
$$0 \rightarrow {\mathbb Z}^{|\partial \overline{S}|} \rightarrow
  \Mod(\overline{S}) \rightarrow \Mod(\widehat{S}) \rightarrow 0.$$
We write
$$\overline{\varepsilon}_{*}\co\overline{G} \rightarrow \widehat{G}$$
to denote the restricted epimorphism.

The inclusion
$$i\co\overline{S} \rightarrow S$$
also induces a homomorphism
$$i_{*}\co\Mod(\overline{S}) \rightarrow \Mod(S).$$
One can show that the kernel of $i_{*}$ is contained in the kernel of
$\varepsilon_{*}$.
We write
$$\overline{i}_{*}\co\overline{G} \rightarrow G$$
to denote the restriction of $i_{*}$ to $\overline{G}$.
By construction, $\overline{i}_{*}$ is surjective.

We also note that
$$\Mod(\widehat{S}) \cong \prod_{r = 1}^{k} \Mod(\widehat{S}_{r})$$
which allows us to view $\widehat{G}$ as a subgroup of the direct product
$$\widehat{G} < \prod_{r = 1}^{k} \widehat{G}_{r}.$$
Denote the projection onto the $r^{th}$ factor by
$$\pi_{r}\co\widehat{G} \rightarrow \widehat{G}_{r}$$
and note that this is surjective.

Suppose now that there exists an isomorphism $\widehat{G}_{r} \rightarrow
{\mathbb F}_{2}$ for some $r$.
$\overline{G}$ is two-generator, hence a quotient of ${\mathbb F}_{2}$,
so we have:
$$\begin{CD}
  {\mathbb F}_{2} @>>> \overline{G} @>{\overline{\epsilon}_{*}}>>
  \widehat{G} @>{\pi_{r}}>> \widehat{G}_{r} @>>> {\mathbb F}_{2}\\
  & & @V{\overline{i}_{*}}VV\\
  & & G\\
  \end{CD}$$
All the arrows are surjections, and free groups are Hopfian (see
\cite{MKS}), so the composition of all the horizontal arrows is an
isomorphism.
Therefore, all horizontal arrows are isomorphisms, and in particular
$\overline{G} \cong {\mathbb F}_{2}$.

Since $\overline{i}_{*}$ is surjective, we'll have $G \cong {\mathbb
F}_{2}$ if $\overline{i}_{*}$ is also injective.  The kernel of
$\overline{i}_{*}$ is contained in the kernel of $\varepsilon_{*}$, and is
therefore contained in the center of $\overline{G}$.  Since $\overline{G}
\cong {\mathbb F}_{2}$, the center is trivial and so $\overline{i}_{*}$
is injective.

Now suppose that $G \cong {\mathbb F}_{2}$, and note that this guarantees
that $\overline{G} \cong {\mathbb F}_{2}$, again appealing to the Hopfian
property.  Because the kernel of $\overline{\varepsilon}_{*}$ is central,
it follows that $\widehat{G} \cong {\mathbb F}_{2}$.  We need to verify
that $\widehat{G}_{r} \cong {\mathbb F}_{2}$ for some $r$.  If this were
not the case, then $K_{r} = \ker (\pi_{r})$ is a non-trivial normal
subgroup of $\widehat{G}$ for each $r$.  An easy induction argument
shows that the commutator group
$$[\ldots[[[K_{1},K_{2}],K_{3}],K_{4}],\ldots,K_{k}]$$
is contained in each $K_{r}$, and hence must be trivial in $\widehat{G}$.
Since $\widehat{G} \cong {\mathbb F}_{2}$, any commutator subgroup of
non-trivial normal subgroups must be non-trivial.
This contradiction proves the proposition.
\end{proof}

\section{Appendix B: Penner's construction} \label{pennerconst}

In \cite{P2}, Penner gives a generalization of a special case of
Thurston's construction for pseudo-Anosov automorphisms.  In this section,
we show that the lower bound $\lambda_{L}$ remains valid for this class of
pseudo-Anosov automorphisms.  In fact, we show that the the dilatations of
pseudo-Anosov automorphisms obtained from this construction are bounded
below by $\sqrt{5} > \lambda_{L}$.

We begin by describing Penner's construction.
Consider $A,B \in {\cal S}'(S)$, label the components $A = a_{1} \cup
\cdots \cup a_{n}$ and $B = b_{1} \cup \cdots \cup b_{m}$, and suppose
that $A \cup B$ fills $S$.
Consider the semi-group ${\frak G}(A,B)$ consisting of all automorphisms
of the form
\begin{equation}
  \label{equ1}
  \prod_{k = 0}^{N} T_{a_{l_{k}}}^{\epsilon_{k}} T_{b_{s_{k}}}^{-\delta_{k}}
\end{equation}
where $N, \epsilon_{k}, \delta_{k} \in {\mathbb Z}_{\geq 0}$.  That is,
${\frak G}(A,B)$ consists of all possible products of positive Dehn twists
about components of $A$ and negative Dehn twists about components of $B$.

There is a subsemigroup ${\frak G}_{0}(A,B)$ consisting of all elements
of ${\frak G}(A,B)$ for which every component of $A$ and $B$ is twisted
along non-trivially at least once in the above product.
In \cite{P2}, Penner proves:

\begin{theorem}[Penner]
\label{pennersthrm}
${\frak G}_{0}(A,B)$ consists entirely of pseudo-Anosov automorphisms.
\end{theorem}

Note that ${\frak G}_{0}(A,B)$ contains all the elements of $\langle
T_{A},T_{B} \rangle$ representable as words in $T_{A}$ and $T_{B}$ where
$T_{A}$ (respectively, $T_{B}$) appears with all positive (respectively,
negative) exponents.  However, this is a relatively small subset of
${\frak G}_{0}(A,B)$, as most elements of ${\frak G}_{0}(A,B)$ do
not obviously lie in $\langle T_{A},T_{B} \rangle$.  Thus, Penner's
construction generalizes a particular case of the construction we have
been considering.

The method which Penner uses to prove Theorem \ref{pennersthrm} allows
one to easily obtain the following bound.

\begin{theorem} \label{pennersdilatation}
The dilatation of any element of ${\frak G}_{0}(A,B)$ is bounded below
by $\sqrt{5}$.
\end{theorem}

The proof we give uses the methods described in \cite{P2}.  We refer the
reader to that paper for a more complete description of those techniques.
We also note that the estimates we give are rough, and this bound is
likely not sharp, though we do not prove this.

\begin{proof}
Fix an element $\phi \in {\frak G}_{0}(A,B)$.  As in \cite{P2}, we
will consider $\phi^{2}$ instead of $\phi$.  Since $\lambda(\phi^{2})
= (\lambda(\phi))^{2}$, it suffices to prove that $\lambda(\phi^{2})
\geq 5$.

At each intersection point of a component $a_{l}$ of $A$ with a component
$b_{s}$ of $B$, apply a homotopy of $b_{s}$ so that it meets $a_{l}$ as in
Figure \ref{pennertrack1}.  The union of the resulting curves is a {\em
bigon track}, $\tau$ (this is essentially a train track except we have
weakened the non-degeneracy condition on complementary regions, allowing
bigons).  Let us denote the branches by $\beta_{1},\ldots,\beta_{K}$.\\

\begin{figure}[ht!]\small
\vglue6pt\cl{\psfig{file=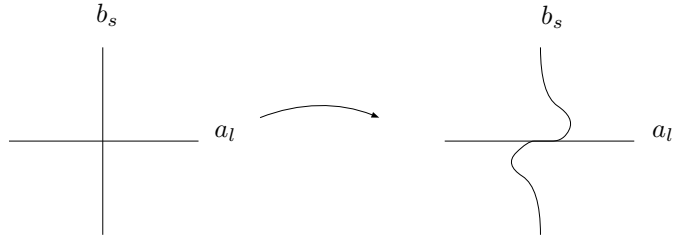,height=1truein}}
\caption{Modifying intersection points}
\label{pennertrack1}
\vglue 2pt
  \setlength{\unitlength}{1in}
  \begin{picture}(0,0)(0,0)
    \put(2.03,1.02){$a_{l}$}
    \put(1.41,1.62){$b_{s}$}
    \put(4.32,1.02){$a_{l}$}
    \put(3.74,1.61){$b_{s}$}
  \end{picture}
\end{figure}

Next, we represent $\phi^{2}$ as a product of Dehn twists in a particular
way so that $\phi^{2}(\tau)$ is easily seen to be carried by $\tau$.
For each component $c$ of $A$ and of $B$ one takes two push-offs,
$c^{\pm}$, one on each side of $c$.  We then express $\phi^{2}$ as a
product of twists along the push-offs, rather than the original curves.
Because every curve which we twist along in $\phi$ shows up twice as
many times in $\phi^{2}$, we can arrange that we twist along {\em both}
push-offs in $\phi^{2}$.  We do this so that we twist along all positive
push-offs in the first application of $\phi$ and then along negative
push-offs in the second application.  Thus, if $\phi$ is given by the
product in (\ref{equ1}), we have
$$\phi^{2} = \prod_{k = 0}^{N} T_{a_{l_{k}}^{-}}^{\epsilon_{k}}
  T_{b_{s_{k}}^{-}}^{-\delta_{k}}\prod_{k = 0}^{N}
  T_{a_{l_{k}}^{+}}^{\epsilon_{k}} T_{b_{s_{k}}^{+}}^{-\delta_{k}}.$$
For each $a_{l}^{+}$, $T_{a_{l}^{+}}(\tau)$ is carried by $\tau$, as is
indicated by Figure \ref{pushtwist1} in the case that $i(a_{l},B) = 1$.
Let us write $M^{a_{l}^{+}}$ to denote the incidence matrix describing
how $\tau$ carries $T_{a_{l}^{+}}(\tau)$.
One can verify that $M^{a_{l}^{+}}$ has the form
$$M^{a_{l}^{+}} = I + R^{a_{l}^{+}}$$
where $I$ is the $K \times K$ identity matrix, and $R^{a_{l}^{+}}$
is a non-negative integral matrix.  Moreover, if $\beta_{p}$ is any
branch contained in $a_{l}$ and $\beta_{q}$ is a branch contained in
$B$ which intersects $a_{l}^{+}$, then the $(p,q)$--entry satisfies
$(R^{a_{l}^{+}})_{pq} = 1$.  Similar statements hold for push-offs
$a_{l}^{-}$, $b_{s}^{+}$, and $b_{s}^{-}$.

\begin{figure}[ht!]\small
\begin{center}
\psfig{file=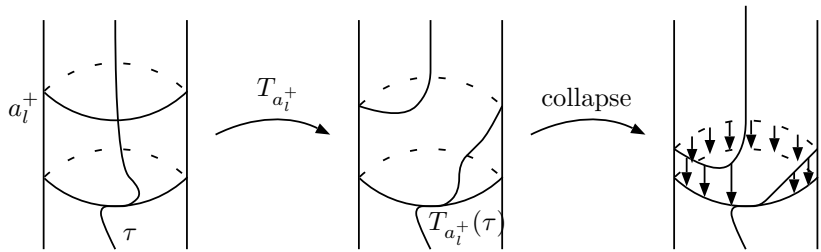,height=1.3truein}
\caption{$\tau$ carrying $T_{a_{l}^{+}}(\tau)$}
\label{pushtwist1}
\end{center}
  \setlength{\unitlength}{1in}
  \begin{picture}(0,0)(0,0)
    \put(.98,.67){\small $\tau$}
    \put(.38,1.32){\small $a_{l}^{+}$}
    \put(1.68,1.42){\small $T_{a_{l}^{+}}$}
    \put(2.58,.72){\small $T_{a_{l}^{+}}(\tau)$}
    \put(3.18,1.37){\small collapse}
  \end{picture}
\end{figure}

In particular, suppose that $a_{l}$ and $b_{s}$
intersect in at least one point $\xi$, and let
$\beta_{i^{+}},\beta_{i^{-}},\beta_{j^{+}},\beta_{j^{-}}$ be the
branches of $\tau$ around $\xi$ as indicated in Figure \ref{pushtwist2}.
The $(p,q)$--entries of $R^{c^{\pm}}$ for $c = a_{l}$ or $b_{s}$ satisfy:
$$\begin{array}{llll}
  (R^{a_{l}^{+}})_{pq} = 1 & \mbox{ for } p = i^{\pm} \, , \, q = j^{+} &
  (R^{a_{l}^{-}})_{pq} = 1 & \mbox{ for } p = i^{\pm} \, , \, q = j^{-}\\
  (R^{b_{s}^{+}})_{pq} = 1 & \mbox{ for } p = j^{\pm} \, , \, q = i^{+} &
  (R^{b_{s}^{-}})_{pq} = 1 & \mbox{ for } p = j^{\pm} \, , \, q = i^{-}\\
\end{array}$$

\begin{figure}[ht!]\small
\begin{center}
\psfig{file=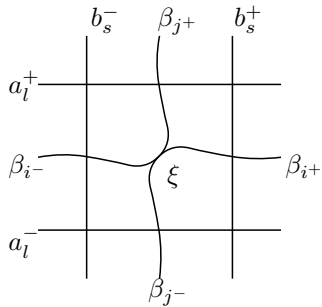,height=1.3truein}
\caption{The branches around the intersection point $\xi$}
\label{pushtwist2}
\end{center}
  \setlength{\unitlength}{1in}
  \begin{picture}(0,0)(0,0)
    \put(3.25,1.17){$\beta_{i^{+}}$}
    \put(1.8,1.17){$\beta_{i^{-}}$}
    \put(2.58,1.94){$\beta_{j^{+}}$}
    \put(2.53,.52){$\beta_{j^{-}}$}
    \put(1.8,1.59){$a_{l}^{+}$}
    \put(1.8,.77){$a_{l}^{-}$}
    \put(2.98,1.94){$b_{s}^{+}$}
    \put(2.23,1.94){$b_{s}^{-}$}
    \put(2.63,1.12){$\xi$}
  \end{picture}
\end{figure}

Now, the incidence matrix $M$ describing $\tau$ carrying $\phi^{2}(\tau)$
is given by the product
$$M = \prod_{k=0}^{N}(M^{a_{l_{k}}^{-}})^{\epsilon_{k}}
  (M^{b_{s_{k}}^{-}})^{-\delta_{k}}
  \prod_{k=0}^{N}(M^{a_{l_{k}}^{+}})^{\epsilon_{k}}
  (M^{b_{s_{k}}^{+}})^{-\delta_{k}}.$$
It is not hard to see that one of $M^{b_{s}^{+}}$ or $M^{b_{s}^{-}}$
occurs between $M^{a_{l}^{-}}$ and $M^{a_{l}^{+}}$ in this product
(these matrices all occur since $\phi \in {\frak G}_{0}(A,B)$).
So, we may write
$$M = X_{0}M^{a_{l}^{-}}X_{1}M^{b_{s}^{\sigma}}X_{2}M^{a_{l}^{+}}X_{3}$$
where $X_{t} = I + Y_{t}$, and $Y_{t}$ is a non-negative integral matrix,
for $0 \leq t \leq 3$, and $\sigma \in \{+,- \}$.
Expanding this out, we see that
\begin{multline*}
M = (I + Y_{0})(I+R^{a_{l}^{-}})(I + Y_{1})(I + R^{b_{s}^{\sigma}})(I +
  Y_{2})(I + R^{a_{l}^{+}})(I + Y_{3}) \\
  = I + R^{a_{l}^{-}} + R^{a_{l}^{+}} + R^{a_{l}^{-}}R^{b_{s}^{\sigma}}
  + R^{a_{l}^{-}}R^{b_{s}^{\sigma}}R^{a_{l}^{+}} + Z
\end{multline*}
where $Z$ is a non-negative integral matrix.
Using the above values for $(R^{c^{\pm}})_{pq}$, one can check that each
of the first $5$ matrices in this last sum has a positive entry in the
$(i^{\pm})$th rows.
It follows that the sum of the entries in each of the $(i^{\pm})$th
rows of $M$ is at least $5$.

The $\beta_{i^{\pm}}$ were arbitrary branches contained in $A$: $a_{l}$
and $b_{s}$ were arbitrary, and every branch in $A$ is adjacent to
some intersection point (eg $(R^{a_{l}^{-}})_{i^{\pm}j^{-}} = 1$,
$(R^{b_{s}^{\sigma}})_{j^{-}i^{\sigma}} = 1$, so the $i^{\pm}i^{\sigma}$
entry of the third term is at least $1$).  Therefore, every row of $M$
with index corresponding to a branch in $A$ has the sum of its entries
being at least $5$.  A similar argument can be made for branches contained
in $B$, and thus it follows that every row of $M$ has sum at least $5$.
Appealing to Theorem \ref{pf1}, we see that the PF eigenvalue of $M$
is at least $5$: take $\vec{U}$ to be the vector with all entries equal
to $1$, and apply the first inequality of the theorem.

The following lemma, which is implicit in the proof of Theorem \ref{pennersthrm} given in \cite{P2} completes the proof.
\end{proof}

\begin{lemma}[Penner]
The PF eigenvalue of $M$ is $\lambda(\phi^{2})$.
\end{lemma}
%
%
%
%

\end{document}